\documentclass{amsart}

\usepackage{amsmath,amsthm,amssymb,amsfonts,eucal,latexsym}
\usepackage{ifthen,epsfig,color,graphicx,rotating,hangcaption,setspace}

\date{\today}

\newtheorem{thm}{Theorem}[section]
\newtheorem{prop}[thm]{Proposition}
\newtheorem{cor}[thm]{Corollary}
\newtheorem{lem}[thm]{Lemma}

\theoremstyle{definition}
\newtheorem{defn}[thm]{Definition}
\newtheorem{exmp}[thm]{Example}

\newtheorem{rem}[thm]{Remark}

\theoremstyle{remark}

\renewcommand{\Re}{\text{Re}}
\renewcommand{\Im}{\text{Im}}
\renewcommand{\(}{\left(}
\renewcommand{\)}{\right)}  
\newcommand{\abs}[1]{\left\lvert #1\right\rvert}
\newcommand{\norm}[1]{\left\| #1\right\|}
\newcommand{\onorm}[1]{\abs{#1}}
\newcommand{\snorm}[1]{\norm{#1}}
\newcommand{\sonorm}[1]{\onorm{#1}}
\newcommand{\vsymb}{\sigma}
\newcommand{\vnorm}[1]{\norm{#1}_{\vsymb}}
\newcommand{\enorm}[1]{\norm{#1}_{e}}

\newcommand{\snbd}[2]{\mathcal{N}({#1},{#2})}
\newcommand{\twovec}[2]{(#1,#2)}
\newcommand{\field}[1]{\mathbb{#1}}
\newcommand{\DD}{\ensuremath{\field{D}}} 
\newcommand{\FF}{\ensuremath{\field{F}}} 
\newcommand{\KK}{\ensuremath{\field{K}}} 
\newcommand{\NN}{\ensuremath{\field{N}}} 
\newcommand{\ZZ}{\ensuremath{\field{Z}}} 
\newcommand{\II}{\ensuremath{\field{I}}}
\newcommand{\IF}{\ensuremath{\field{I}\field{F}}} 
\newcommand{\IR}{\ensuremath{\field{I}\field{R}}} 
\newcommand{\CC}{\ensuremath{\field{C}}} 
\newcommand{\Ct}{\ensuremath{\field{C}^{{2}}}}
\newcommand{\Cn}{\ensuremath{\field{C}^{{n}}}}
\newcommand{\RR}{\ensuremath{\field{R}}}
\newcommand{\Rt}{\ensuremath{\field{R}^{{2}}}}
\newcommand{\Rn}{\ensuremath{\field{R}^{{n}}}}
\newcommand{\Rtn}{\ensuremath{\field{R}^{{2n}}}}

\newcommand{\eps}{\varepsilon}
\newcommand{\Henon}{H\'{e}non }
\newcommand{\Hypatia}{Hypatia}

\newcommand{\boxch}{box chain }
\newcommand{\boxchcn}{box chain construction }
\newcommand{\boxcov}{\mathcal{B}} 
\newcommand{\boxchrec}{box chain recurrent } 
\newcommand{\boxchrecset}{box chain recurrent set } 
\newcommand{\boxchrecsets}{box chain recurrent sets } 
\newcommand{\boxchmod}{box chain model } 
\newcommand{\boxchmods}{box chain models } 
\newcommand{\boxchrecmod}{box chain recurrent model } 
\newcommand{\boxchrecmods}{box chain recurrent models }
\newcommand{\boxchtr}{box chain transitive }
\newcommand{\boxchtrcomp}{box chain transitive component }
\newcommand{\boxchtrcomps}{box chain transitive components }

\newcommand{\chainrec}{\mathcal{R}} 
\newcommand{\chrec}{\chainrec}
\newcommand{\chrecdeltapn}{\chainrec({{\delta}'_n})}
\newcommand{\chrecdeltapnmo}{\chainrec({{\delta}'_{n-1}})}
\newcommand{\chrecdeltapo}{\chainrec({{\delta}'_0})}
\newcommand{\chrecdeltap}{\chainrec({{\delta}'})}
\newcommand{\chrecpdeltap}{\chainrec'({{\delta}'})}
\newcommand{\chrecdelta}{\chainrec({{\delta}})}
\newcommand{\chrecepspn}{\chainrec({{\eps}'_n})}

\newcommand{\chreceps}{\chainrec({\eps})}
\newcommand{\chrecepsp}{\chainrec({\eps'})}
\newcommand{\chrecpepsp}{\chainrec'({{\eps'}})}
\newcommand{\chreceta}{\chainrec({{\eta}})}
\newcommand{\chrecnu}{\chainrec({{\nu}})}

\newcommand{\inKsel}{sink basin }
\newcommand{\inKselsub}{sink basin subdivision }

\newcommand{\showcomments}{yes}

\newsavebox{\commentbox}
{\ifthenelse{\equal{\showcomments}{yes}}%
{\footnotemark
    \begin{lrbox}{\commentbox}
    \begin{minipage}[t]{1in}\raggedright\sffamily\small
    \footnotemark[\arabic{footnote}]}
{\begin{lrbox}{\commentbox}}}%
{\ifthenelse{\equal{\showcomments}{yes}}%
{\end{minipage}\end{lrbox}\marginpar{\usebox{\commentbox}}}
{\end{lrbox}}}

\newcommand{\drawfigaltpertwoA}{\scalebox{.75}{\includegraphics{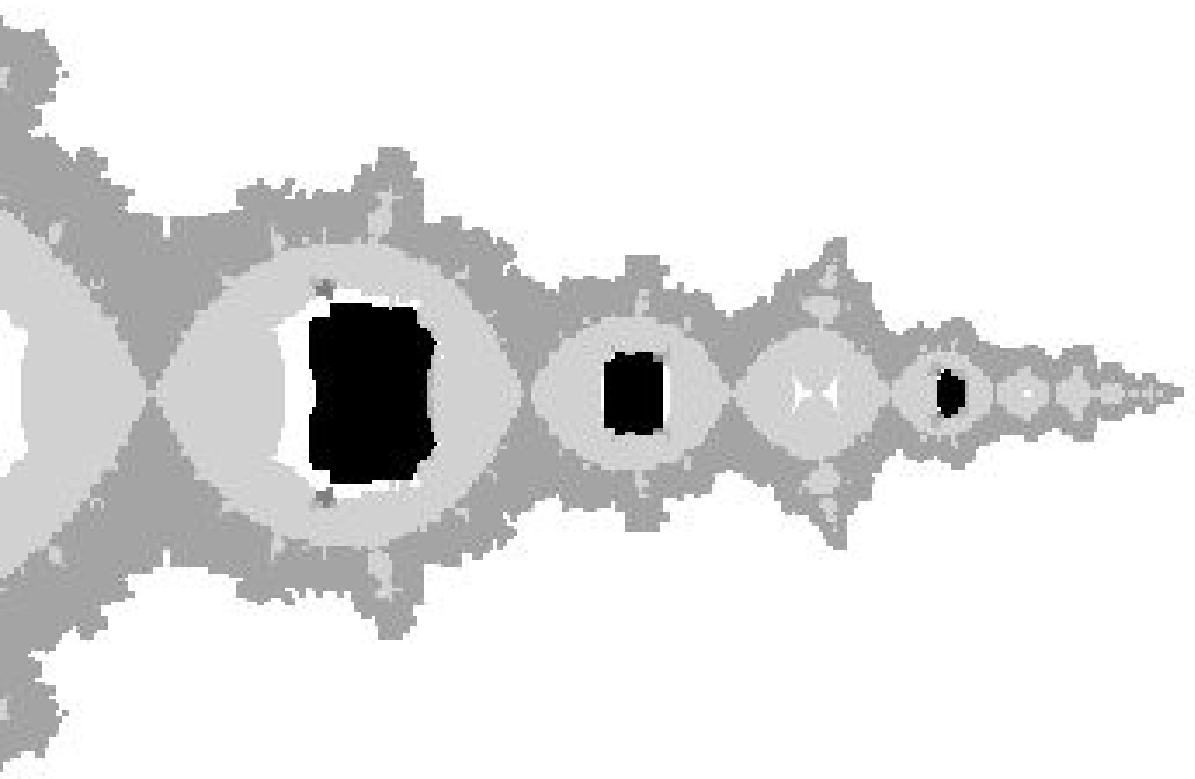}}}
\newcommand{\drawfigaltpertwoC}{\scalebox{.75}{\includegraphics{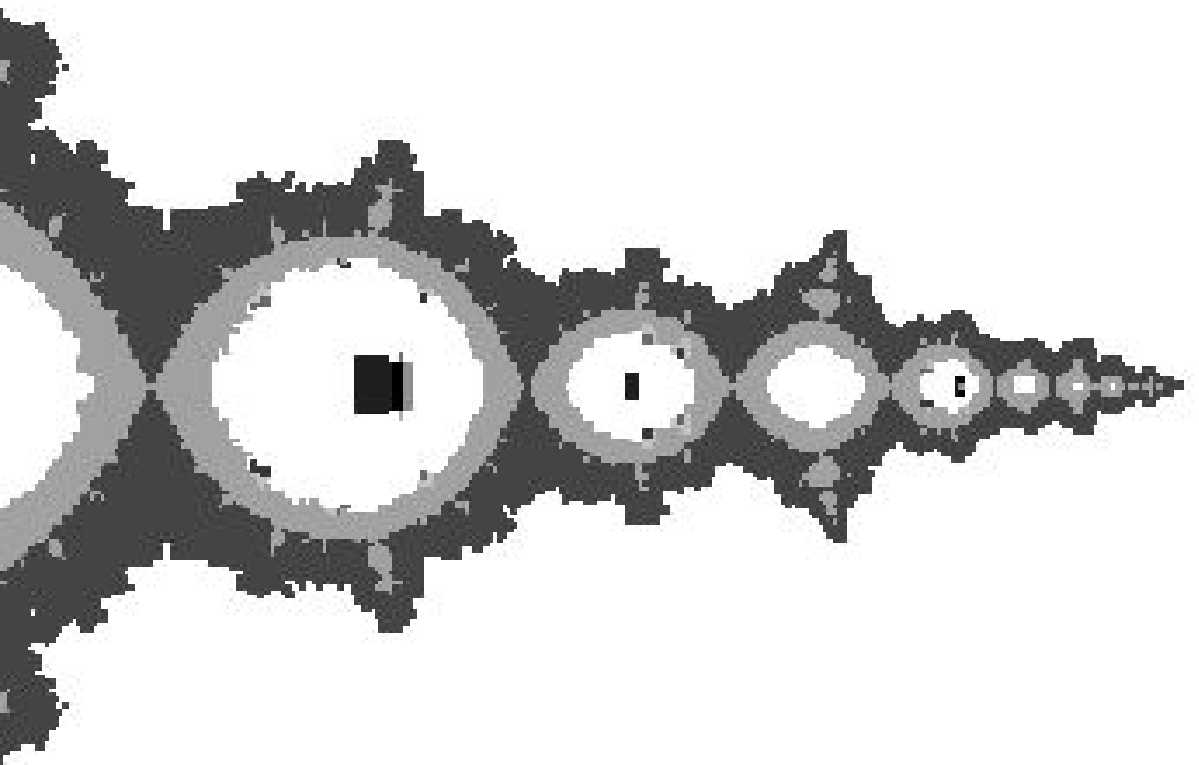}}}

\newcommand{\drawfigric}{\scalebox{.25}{\includegraphics{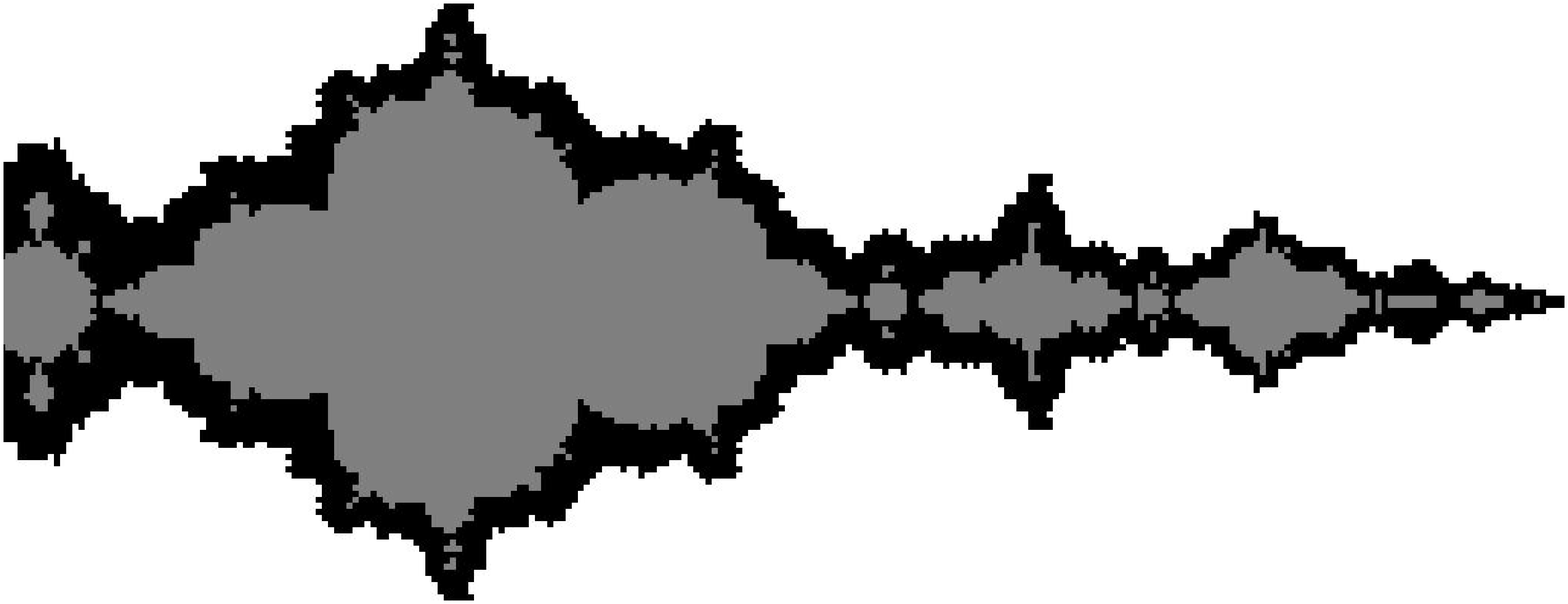}}}

\newcommand{\drawfigcubdoublebas}{\scalebox{.7}{\includegraphics{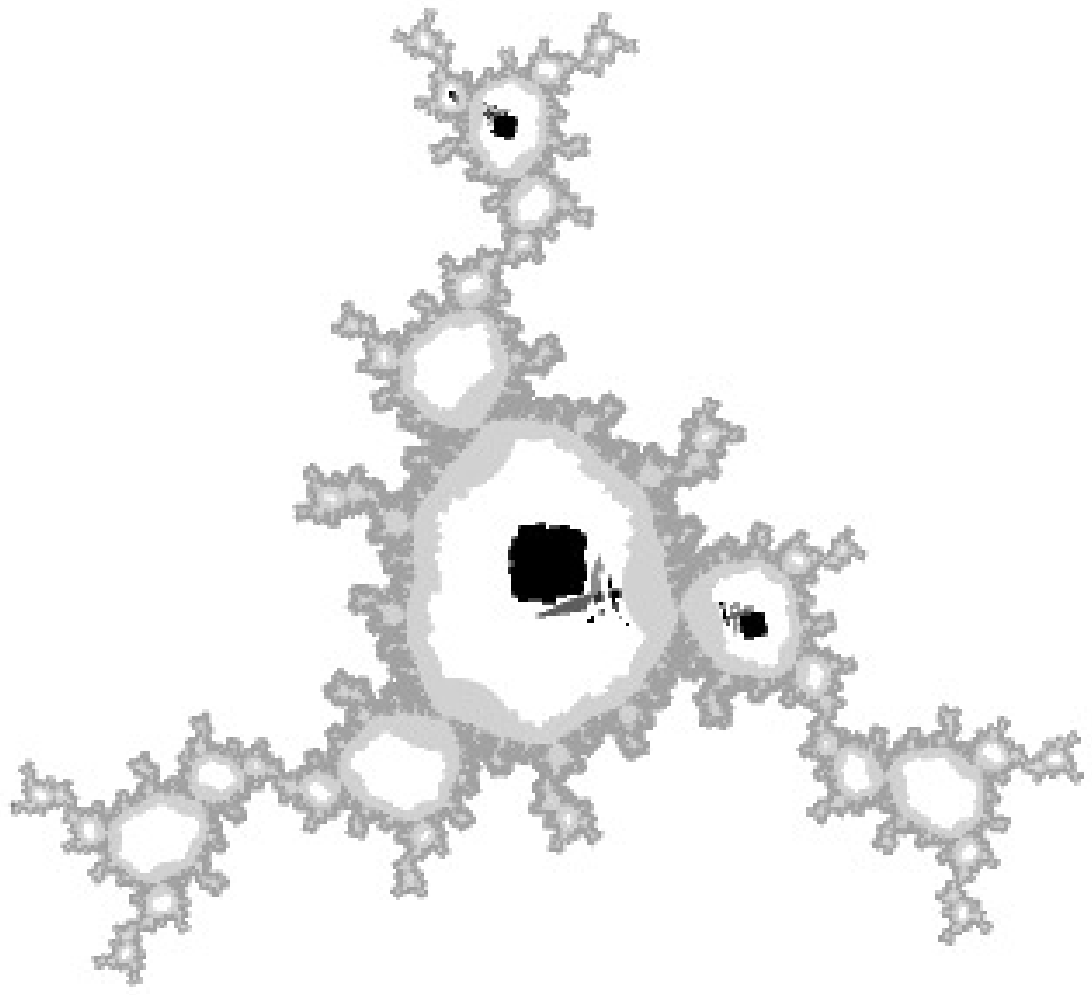}}}

\newcommand{\drawfigcomplexhorseFA}{\scalebox{.25}{\includegraphics{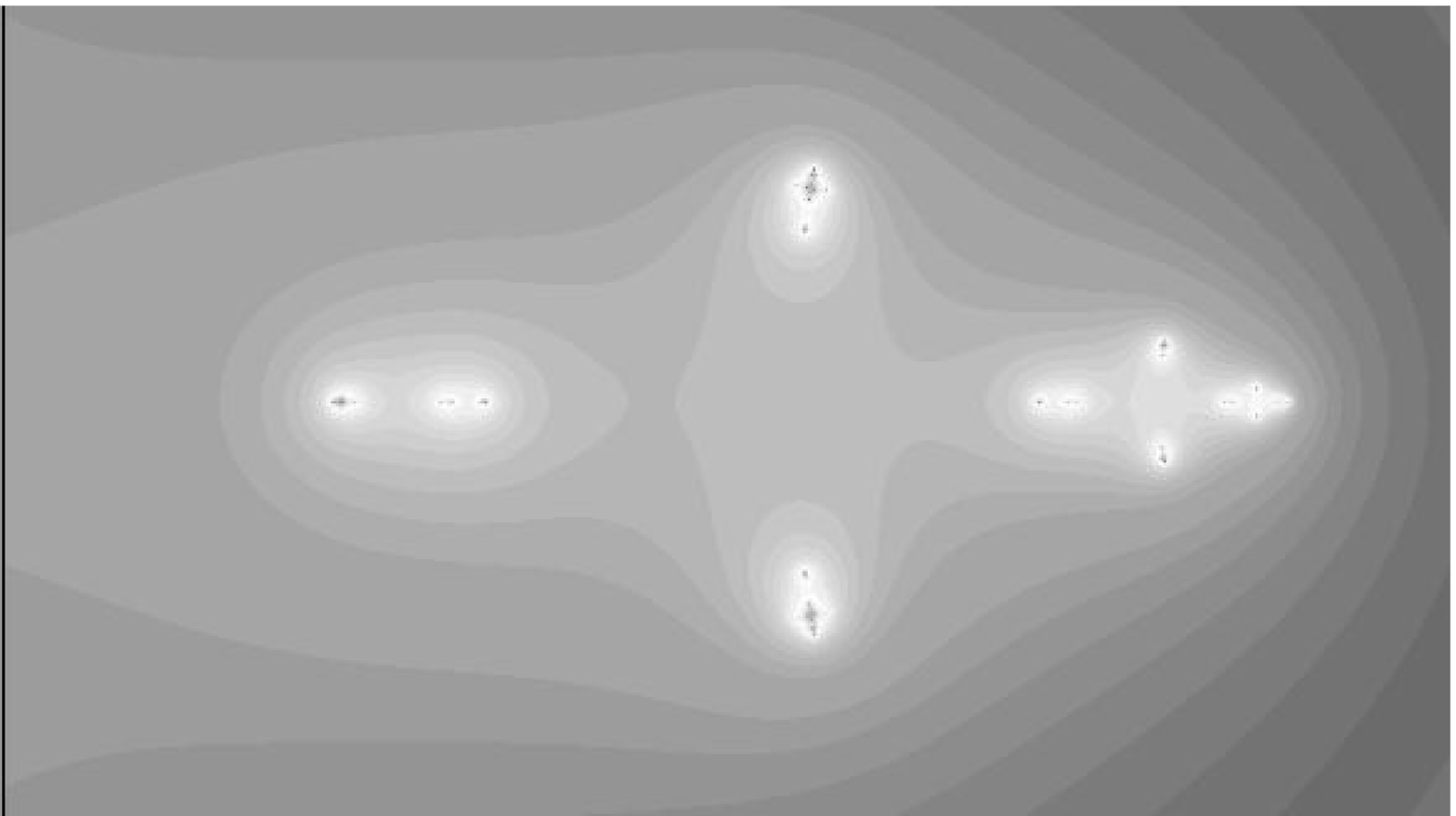}}}
\newcommand{\drawfigcomplexhorse}{\scalebox{.5}{\includegraphics{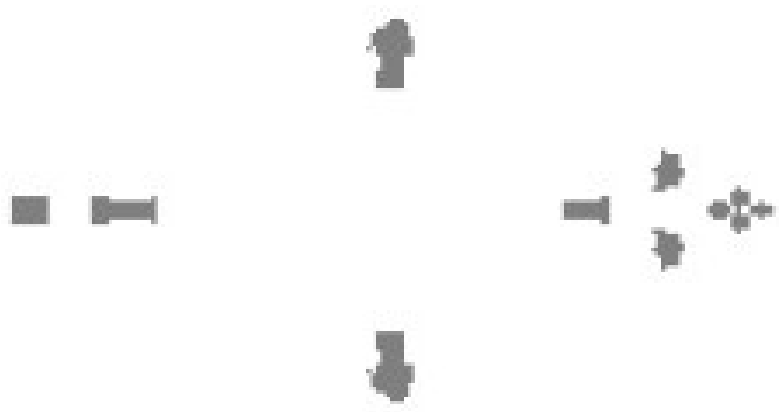}}}

\newcommand{\drawfigHrealhorse}{\scalebox{.23}{\includegraphics{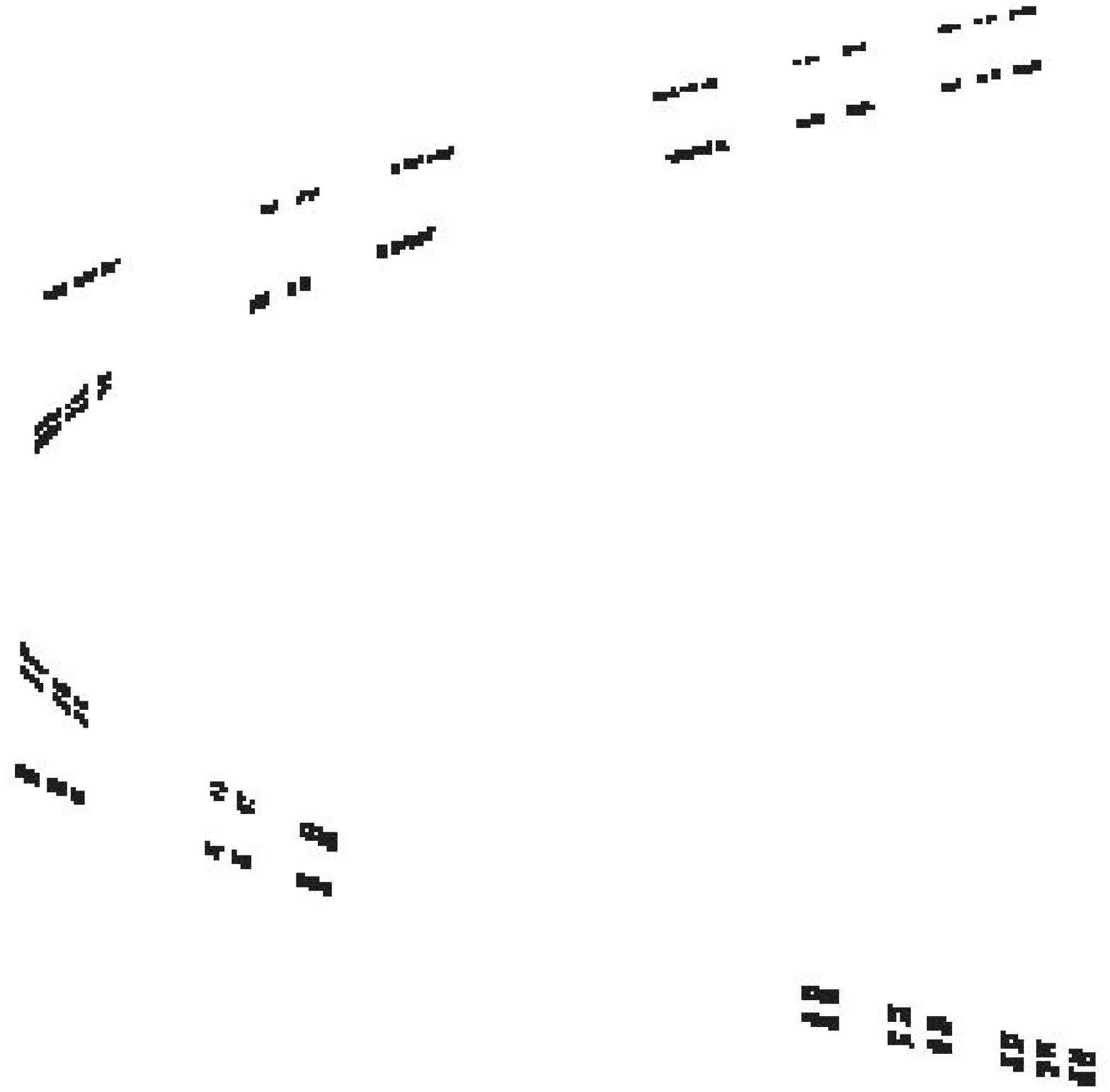}}}
\newcommand{\drawfigHrealhorseW}{\scalebox{.3}{\includegraphics{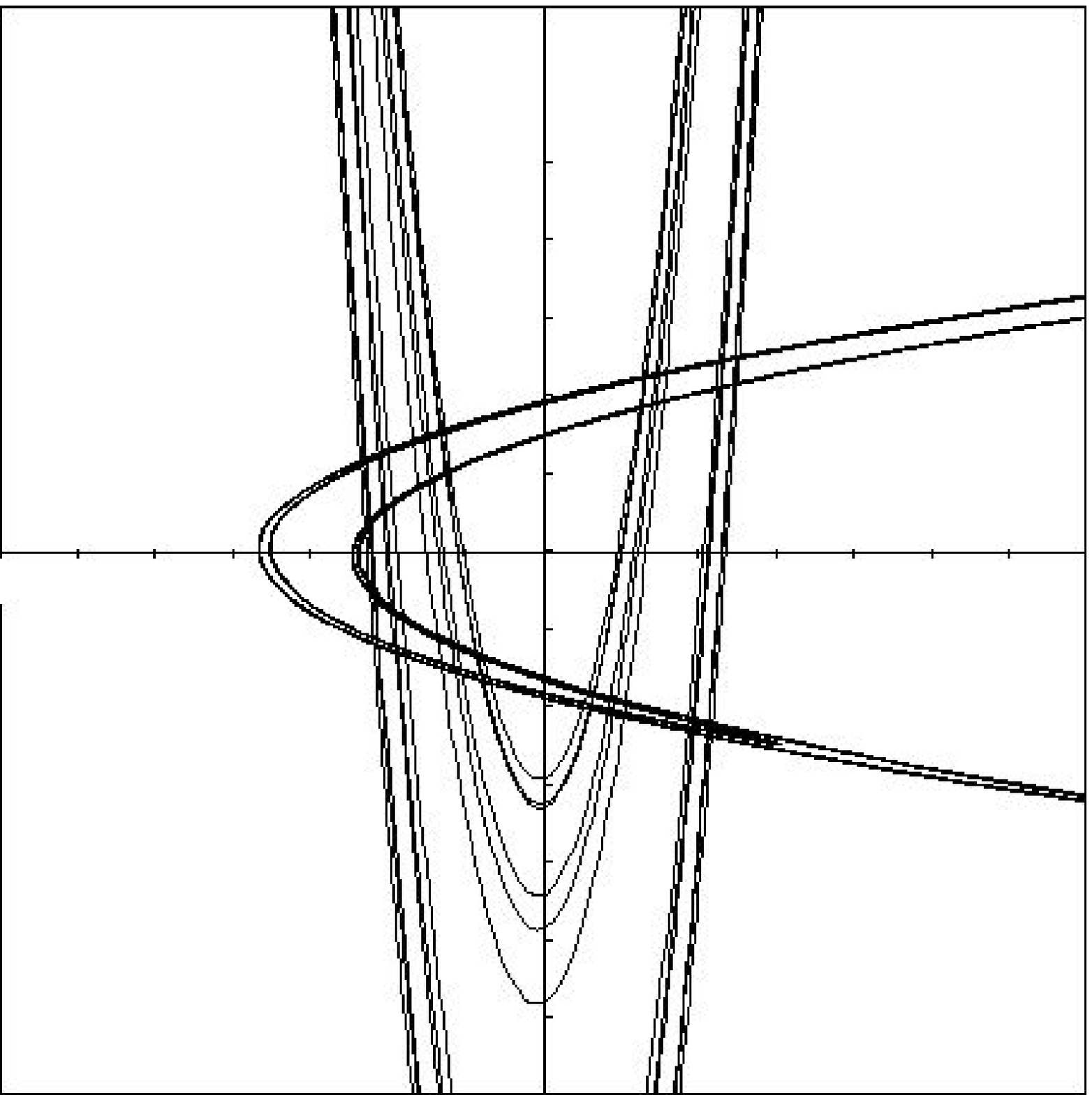}}}

\title[Modeling chain recurrence for H\'{e}non maps]
{Rigorous numerical models \\ for the dynamics of
 complex H\'{e}non mappings \\ on their chain recurrent sets
}

\author[S.L. ~Hruska]{}

\email{shruska@msm.umr.edu}



\begin{document}

\maketitle

\centerline{\scshape  Suzanne Lynch Hruska}
\medskip

{\footnotesize
\centerline{ Department of Mathematics }
\centerline{ Indiana University }
\centerline{ Rawles Hall }
\centerline{ Bloomington, IN 47405, USA }
}
\medskip


\bigskip

\begin{quote}{\normalfont\fontsize{8}{10}\selectfont
{\bfseries Abstract.}
We describe a rigorous and efficient computer algorithm for  building 
a model of the dynamics of a polynomial diffeomorphism of 
$\Ct$ on its chain recurrent set, $\chainrec$, and for sorting 
points into approximate chain transitive components.  Further, we give explicit 
estimates which quantify how well this algorithm approximates the chain 
recurrent set and distinguishes the chain transitive components. We also 
discuss our implementation for the family of  \Henon mappings, 
$f_{a,c}(x,y) =  (x^2 + c - ay, x)$, 
into 
a computer program called \textit{Hypatia}, and give several examples 
of running Hypatia on \Henon mappings.
\par}
\end{quote}

\maketitle


\section{Introduction}
\label{sec:intro}

Computer work, especially computer graphics, has been an important
tool of discovery in the field of dynamical systems. This paper is also
concerned with the use of computers but it has a different goal. The
goal of this paper is to rigorously establish some results on the location of the set
of points where recurrent behavior takes place.
We focus on the class of polynomial diffeomorphisms of $\Ct$, which includes
the widely studied family of \Henon mappings, $f_{a,c}(x,y) = (x^2+c-ay, x)$.

We start by examining a
rigorous and effective computer algorithm for 
building a model of the dynamics of a map, $f$, on its chain recurrent set, $\chainrec$, 
and for sorting points into approximate chain transitive components.
We call this algorithm the \textit{\boxch construction}.
%
The same basic algorithm has been 
studied previously, in different settings.  Osipenko and Campbell (\cite{Osi, OsiCamp}) approximate the chain recurrent set for a homeomorphism of a smooth, real, compact manifold. Eidenschink (\cite{Eiden}) discusses a similar procedure for real flows.  A philosophically related procedure is studied in \cite{DellHoh, Dell1}, but their case of interest is the attractor of a real map (rather than the chain recurrent set).  \cite{KMisch} is a recent survey of this work. 
The focus in the previous body of work is to develop a very general  procedure for rigorously approximating $\chrec$ for continuous maps or flows in $\Rn$.  

In this paper, we restrict our attention to polynomial diffeomorphisms of $\Ct$, which allows us to adapt the \boxchcn and its implementation to be more efficient.  In addition, we establish estimates on the accuracy of our model.  These are explicit estimates, involving only the inputs to the program, which quantify how closely our \boxchcn approximates the chain recurrent set.  This allows us to predict when an execution of the program will successfully separate the distinct chain transitive components.
We contrast this with the work of
 Dellnitz and Hohmann (\cite{DellHoh}), which gives results on the accuracy of the approximation in the case that the map is hyperbolic; however, their estimates involve constants of hyperbolicity, which they do not discuss how to calculate.  Osipenko and Campbell (\cite{OsiCamp}) give approximation estimates, in terms of some of the inputs to the program as well as selected output, hence the accuracy of their approximation can only be measured after execution of the algorithm.  
%

%
%
%
Our first result on accuracy of the box chain construction for polynomial diffeomorphisms of $\Ct$ can be paraphrased as follows.

\begin{thm} \label{thm:babyRinBoxCov}
Let $f$ be a polynomial diffeomorphism of $\Ct$ (of dynamical degree $d > 1$). Suppose  $\boxcov_0$ is a closed box in $\Ct$ containing the $\delta'_0$-chain recurrent set, $\chrecdeltapo$, for some $\delta'_0 >0$ (for example,  take $\boxcov_0,$ and $\delta=\delta'_0 >0$ as in Proposition~\ref{prop:RinV}). 

The \boxchcn produces sequences of constants, $\{ \delta_n \}$ and $\{ \eps_n \}$,  directed graphs, $\{ \Gamma_n \}$, and regions in $\Ct$, $\{ \boxcov_n \}$, for $n \geq 1$, such that
\begin{enumerate}
\item $\delta_n \ll \eps_n$ and both $\eps_n \downarrow 0$ and $\delta_n \downarrow 0$ as $n\to \infty$,
\item the vertex set of $\Gamma_n$ is a collection of closed boxes $\{ B^n_k \}$ in $\Ct$, which have side length at most~$\eps_n$,
\item  each $\boxcov_n$ is the region in $\Ct$ defined by the union of  the $B^n_k$,
\item  $\{ \boxcov_n \}$  is nested, \textit{i.e.},
$ \{ \cdots \subset \boxcov_n \subset \cdots \subset \boxcov_1 \subset \boxcov_0 \}$,
\item there is guaranteed to be an edge in $\Gamma_n$ from $B^n_k$ to $B^n_j$ if $f(B^n_k)$ intersects a $\delta_n$-neighborhood of $B^n_j$,  
and 
\item we can calculate 
explicit sequences $\{ \eps'_n \}_{n=1}^{\infty}$ and $\{ \delta'_n \}_{n=1}^{\infty}$, and an explicit constant $C$, such that for  every $n\geq 1$,
\begin{enumerate}
\item $\eps'_n \downarrow 0$, in particular $\eps_n < \eps'_n \leq \delta_n + C \eps_n$,
\item  $\delta'_n$ is nonincreasing and converges to zero, with $\delta'_n < \delta_n$, 
and
\item
$\chrecdeltapn  \subset \boxcov_n \subset \chrecepspn.$
 \end{enumerate}
\end{enumerate}
\end{thm}

\begin{defn} \label{defn:boxterms}
For any $n \geq 1$, suppose $(\eps, \delta, \Gamma, \boxcov) = (\eps_n, \delta_n, \Gamma_n, \boxcov_n)$ are produced by the \boxchcn at step $n$, and satisfy Theorem~\ref{thm:babyRinBoxCov}.

We call the region $\boxcov$  an {\em $(\eps, \delta)$-\boxchrec set},
and the graph $\Gamma$ an {\em $(\eps, \delta)$-\boxchrecmod of $f$}.
Each edge-connected component $\Gamma'$ of $\Gamma$ is called an {\em $(\eps, \delta)$-\boxchtr component.}
\end{defn}

Conclusions (1) through (5) of Theorem~\ref{thm:babyRinBoxCov} follow immediately from the description of the \boxch construction, given in Section~\ref{sec:construction}.  Conclusion (6) is our first significant a priori estimate on the accuracy of our model, and is established in Sections~\ref{sec:trap} through~\ref{sec:resboxcov}.  There we show that for the case of \Henon mappings, for any chosen $  
R' > \left[ 1+\abs{a} + \sqrt{(1+\abs{a})^2 +4\abs{c}} \right]/2,$
we get 
\begin{eqnarray*}
\delta'_0 & = & \left[ (R')^2 - (1+\abs{a})R' - \abs{c} \right]/2,\\
\delta'_n & = & \min \left(\delta'_0, \ \left[-(2R'+\abs{a}+1) +\sqrt{(2R'+\abs{a}+1)^2+4\delta_n} \right] \right), and\\
\eps'_n & = & \delta_n + \eps_n(1 + \abs{a} + 2R') + \eps_n^2. 
\end{eqnarray*}

A {\boxchrecmod } of $f$ satisfies the definition of a {\em symbolic image of $f$}, given by Osipenko in \cite{Osiold}.
Osipenko and Campbell (\cite{OsiCamp}) derive similar estimates to (6), but their version of $\eps'_n$ depends on measuring the size of the images of boxes $B_k$, after they have been computed.  

We have implemented our efficient \boxchcn for \Henon mappings into a computer program we call Hypatia. Using Hypatia we have produced many examples of 
\boxchrecmods of \Henon mappings.
\begin{footnote}
{Write to the author to obtain a copy of this C++, unix program.}
\end{footnote}
The \boxchcn also has an immediate analog for polynomial maps of $\CC$, which we include in Hypatia for quadratic and cubic polynomials.
We keep the arithmetical computations in Hypatia rigorous using  \textit{interval arithmetic with directed rounding}. This method was recommended to us by Warwick Tucker, who used it in his recent computer proof that the Lorenz differential equation has the conjectured geometry (\cite{War}).  
See Appendix~\ref{sec:IA} for a brief introduction to interval arithmetic.

The examples we are most interested in studying with Hypatia are \Henon mappings which are \textit{hyperbolic}.  These are the class of maps which have the ``simplest'' chaotic dynamics, and are in fact stable under small perturbation. Thus hyperbolic maps are the most amenable to rigorous computer investigation.  In fact, Bedford and Smillie (\cite{BS2}) have shown that for hyperbolic polynomial diffeomorphisms of $\Ct$, the chain recurrent set
 is well-behaved, in that it consists of simply the Julia set together with finitely many attracting or repelling periodic points.

The simplest hyperbolic \Henon mappings can be described in terms of the dynamics of some quadratic polynomial.  In fact, if $f_{a,c}$ is a \Henon mapping with
$a$ sufficiently small and $c$ is such that the polynomial $P_c(z) = z^2+c$ is hyperbolic, then $f|_J$ is topologically conjugate to the
function induced by $P$ on the inverse limit
${\lim_{\leftarrow} (J, P)}$ (\cite{HOV2}).  In this case, we say that  {\em $f$ is described by $P$}, or simply that {\em $f$ exhibits one dimensional behavior}.  
The work of Hubbard and Papadontonakis  (\cite{HubKarl, CUweb}), and more specifically the work of Oliva (\cite{OT}), suggests points in
parameter space which are conjectured to be hyperbolic and to have
interesting properties. At the moment though neither the hyperbolicity nor
the interesting properties have been established rigorously.

Thus a significant problem in the study of the \Henon family is to understand which maps exhibit one dimensional behavior, and to describe the behavior of maps which do not.
Motivated by this problem, in this paper we use the \boxchcn and its implementation in Hypatia to 
build \boxchrecmods for several interesting examples of \Henon mappings.
Further, we use the results of this paper as the first step in a study of the property of hyperbolicity  for polynomial diffeomorphisms of $\Ct$ in \cite{SLHthree}, and in a study of hyperbolicity (\textit{i.e.}, expansion) for polynomial maps of $\CC$ in \cite{SLHtwo}.   \cite{SLHT} contains an earlier version of the work of this paper, as well as that of \cite{SLHthree, SLHtwo}.

\begin{exmp} \label{exmp:introaltpertwo}
One of the simplest \Henon mappings
which appears not to exhibit one dimensional behavior 
is $f_{a,c}$ with $(a,c) = (.125, -1.24)$.   
Oliva (\cite{OT}) gave combinatorial evidence that this diffeomorphism is hyperbolic with a period two 
attracting cycle, but is not described by a quadratic polynomial.  Using our program Hypatia, we applied the \boxchcn to the nearby map $(a,c) = (.15, -1.1875)$, which seems to be topologically conjugate to $(a,c) = (.125, -1.24)$.  We computed a sequence terminating at a \boxchrecset $\boxcov_{8}$. 
%

For a qualitative estimation of the accuracy of a \boxchrec set, $\boxcov$, we can sketch the intersection of $\boxcov$ with a dynamically significant one-dimensional submanifold of $\Ct$: the unstable manifold of a saddle fixed point, which has a natural parameterization by $\CC$.  This process is explained in Appendix~\ref{sec:henwuppics} and Section~\ref{sec:drawboxcovwup}.  
In Figure~\ref{fig:altpertwo} 
we use this procedure to illustrate two \boxchrec sets, $\boxcov_{7}$ and $\boxcov_{8}$, for the \Henon mapping with  $(a,c) = (.15, -1.1875)$.

\begin{figure}
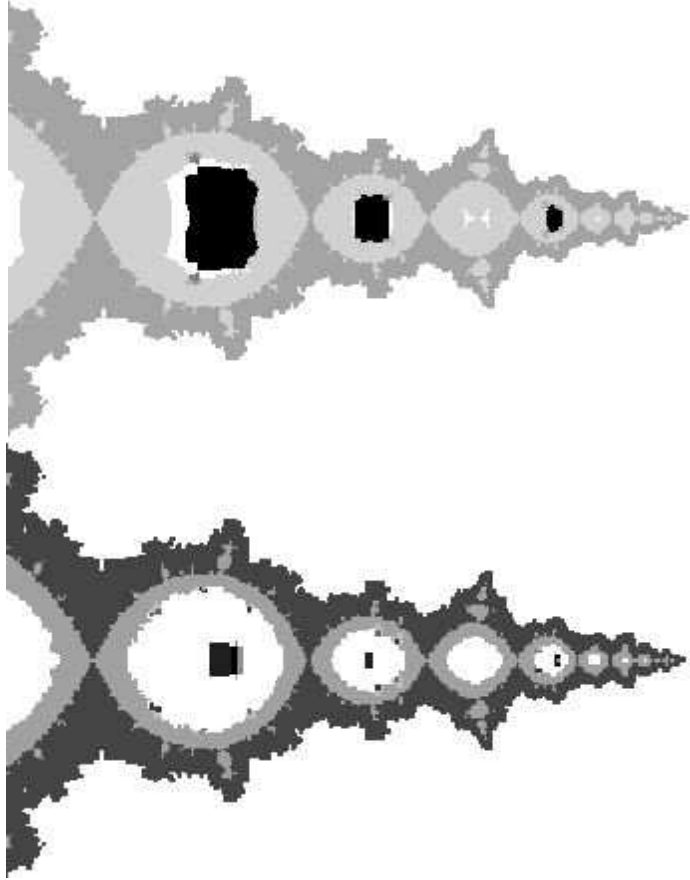

\begin{center}
\drawfigaltpertwoA
\ \ \ \drawfigaltpertwoC
\caption{For $f_{a,c}$ the \Henon mapping with  $c=-1.1875, a=.15$, $\chainrec$ appears to be $J$ and a period $2$ sink.  Shown in this figure are two \boxchrec sets, restricted to the unstable manifold of a saddle fixed point, with its natural parameterization. 
On the top,  boxes are of side length $2R/2^6$ and $2R/2^7$, where $R = 1.9$.
On the bottom is a refinement obtained from subdividing once some of the boxes on the top.  
The innermost, darkest regions in these figures are contained in the \boxchtrcomp containing the 2-cycle.  Each \boxchtrcomp containing $J$ is shaded two tones by a heuristic algorithm, in order to illustrate how close the component is to $J$. 
The top is the crudest \boxchrecmod which separates $J$ from the sink.  
In the bottom figure, the small dark spots skirting the inner circles of the neighborhood of $J$ form a \boxchtrcomp of points which are pseudorecurrent but not recurrent.
} 
\label{fig:altpertwo} 
\end{center}
\end{figure}
%

\end{exmp}

Assuming the conjectural dynamics holds, the \boxchrecmods constructed for Example~\ref{exmp:introaltpertwo} are both successful in the sense of the following.

\begin{defn}  \label{defn:sep}
Let $f$ be a polynomial diffeomorphism of $\Ct$ (of dynamical degree $d > 1$).
Let $\Gamma$ be a \boxchrecmod of $f$.  We call $\Gamma$ {\em separating} if there are two chain transitive components, $\chainrec^j$ and $\chainrec^k$ (of $\chrec$), which lie in different \boxchtrcomps of $\Gamma$.   In this case we say $\Gamma$ {\em separates } $\chainrec^j$ and $\chainrec^k$. 

Further, we call $\Gamma$ {\em fully separating} if it separates every pair of chain transitive components.

 %
\end{defn}

%

\begin{exmp} \label{exmp:introper31}
Another interesting example studied by Oliva (\cite{OT}) is the \Henon mapping $f_{a,c}$ with $(a,c) = (.3, -1.17)$.  
%
%
We call this the \textit{3-1-map}, because it appears to be hyperbolic, with $\chainrec$  
 consisting of three chain transitive components: $J$, an attracting fixed point, and an attracting cycle of period  three.  In contrast, quadratic polynomials cannot have more than one attracting cycle, thus this map appears not to exhibit one dimensional behavior.
We applied the \boxchcn to the 3-1-map, but were surprised to be unable to obtain a separating \boxchrec model.  The best \boxchrecset we obtained is shown in Figure~\ref{fig:per31}, intersected with the parameterized unstable manifold of a saddle fixed point.
\begin{figure}
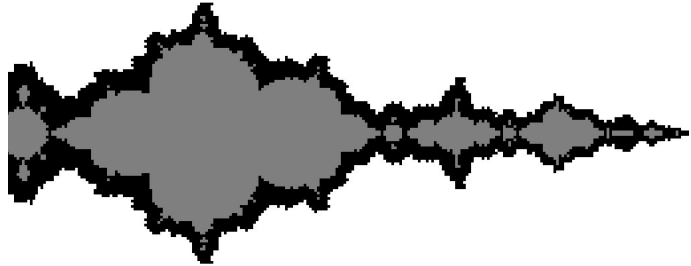

\begin{center}
\drawfigric
\caption
{A \boxchrecset for the \Henon mapping $f_{a,c},$ with  $a=.3, c=-1.17$, restricted to the unstable manifold of a saddle fixed point, with its natural parameterization. Here boxes are of 
side length between $2R/2^7$ and $2R/2^8$, where $R=2.01$.  Lighter gray points were heuristically found to be in $K^+$. Unfortunately, this \boxchrecmod does not separate $J$ from either the fixed sink or the attracting three-cycle.}
\label{fig:per31}
\end{center}
\end{figure}

\end{exmp}

Our difficulties with the 3-1-map motivated the following theorem, in which we calculate explicit bounds on $\eps$ and $\delta$ to guarantee that an $(\eps, \delta)$-\boxchrecmod will separate the fixed sink from the $3$-cycle and the Julia set. This gives a theoretical quantification of the computational difficulty of studying the 3-1-map. 

\begin{thm} \label{thm:babysinkJsep}
Suppose $f_{a,c}$ is a \Henon mapping with an attracting fixed point $p$, with $\lambda_1 \neq \lambda_2$ eigenvalues of $D_pf$, and $\lambda = \max (\abs{\lambda_1}, \abs{\lambda_2})$. 
Set
\[
\tau = \frac{ \abs{\lambda_1 -\lambda_2}^2}
{(2+\abs{\lambda_1}+\abs{\lambda_2})(2+\lambda^2+\abs{a})}.
\]

Let $\Gamma$ be an $(\eps, \delta$)-\boxchrecmod of $f$. 
Let $M>1$ satisfy $\delta < \eps/M$.
Set
\[
\kappa 
=
\left[ 1 + 1/M + \max \{1, (1-\lambda)\sqrt{\tau} +2\snorm{p} + \sonorm{a} 
\} \right]. 
\]

If $\eps < \frac{1}{2} \( -\kappa + \sqrt{\kappa^2+\tau(1-\lambda)^2 } \),$
then $\Gamma$ separates the fixed sink from every other chain transitive component of $\chainrec$.
\end{thm}

This  theorem applied to the 3-1-map yields that boxes of side length less than $4 \times 10^{-5}$ would guarantee separation.  However this is several orders of magnitude smaller than current resources allowed us to achieve with Hypatia.  
This demonstrates the need for the development of a more sophisticated construction for rigorously approximating chain recurrent sets of complex \Henon mappings.
 
Finally, we outline the remaining sections.   In Section~\ref{sec:boxcover} we describe the \boxchcn for polynomial diffeomorphisms of $\Ct$, and prove Theorem~\ref{thm:babyRinBoxCov} by calculating explicit estimates on how well a \boxchrecset approximates the chain recurrent set. 
In Section~\ref{sec:codetricksone} we use some dynamical information about the map to develop two enhancements to the basic construction, significantly improving computational efficiency.
In Section~\ref{sec:sinkJsep} we show a theoretical limitation of the \boxchcn by proving Theorem~\ref{thm:babysinkJsep}, and applying the estimates of the theorem to the 3-1-map. In Section~\ref{sec:examples} we discuss examples generated with Hypatia, for \Henon mappings and a polynomial map of $\CC$.  We have included relevent background material on  the chain recurrent set and the dynamics of \Henon mappings in Appendix~\ref{sec:bckgnd}.  In Appendix~\ref{sec:IA} we sketch the basics of interval arithmetic.

\medskip

\noindent
\textbf{Acknowledgements.}
We thank John Smillie for providing guidance on this project,  John Hubbard, Greg Buzzard, and Warwick Tucker for many helpful conversations on the topic,
James Yorke, John Milnor, and Eric Bedford for advice on the preparation of this paper, Robert Terrell for technical support, and Michael Benedicks for pointing out to us that \cite{KMisch} describes a procedure similar to ours.

\section{The \boxchcn for polynomial diffeomorphisms of $\Ct$}
\label{sec:boxcover}

In this section we start by giving an outline of the \boxch construction, then discuss how we carry it out for polynomial diffeomorphisms of $\Ct$, to calculate the estimates leading to Theorem~\ref{thm:babyRinBoxCov}.  To calculate our estimates, we assume a polynomial diffeomorphism of $\Ct$ (of dynamical degree $d(f)>1$) is a finite composition of {\em generalized \Henon mappings}, which are maps of the form $f(x,y) = (p(x)-ay,x)$, for $p$ monic of degree greater than one (see Appendix~\ref{sec:bckgnd}).


\subsection{Efficient neighborhoods}
\label{sec:metric}

Before we begin our theoretical calculations, we want to specify that we do not use the euclidean metric. 
It is more natural for  computer calculations 
to consider vectors in $\Rtn$,
rather than $\Cn$, and use the $L^{\infty}$ metric, rather than 
euclidean.
Thus throughout this paper, $\snorm{\cdot}$ will denote the $L^{\infty}$ 
norm on~$\Rtn$, so that for a vector $x=(x_1,\ldots,x_n)\in\Cn$,
\[
\snorm{x}=\max\{ \abs{\Re(x_k)},\abs{\Im(x_k)} \colon 1\leq k\leq n\}.
\]
Also, let ${\snbd{S}{r}}$ denote the open $r$-neighborhood about the set
$S$ in the metric $d_{\infty}$ induced by $\snorm{\cdot}$, \textit{e.g.},
\[
{\snbd{0}{r}} = \{x=(x_1,\ldots,x_n)\in\Cn \colon
	\abs{\Re(x_k)}<r \text{ and } \abs{\Im(x_k)}< r \}.
\]
We use the simpler notation $\sonorm{\cdot}$ for dimension $n=1$.
This metric is {\em uniformly
equivalent} to the euclidean metric on~$\Cn, \enorm{\cdot}$, since
$
\frac{1}{\sqrt{2n}} \enorm{x} \leq \snorm{x} \leq \enorm{x}.
$
Neighborhoods are slightly different with respect to two uniformly equivalent
norms, but the topology generated by them is exactly the same, thus they
can practically be used interchangeably.
Similarly, the $\eps$-chain recurrent set $\chreceps$ depends 
on choice of metric, but since $\chainrec^{e}({\eps}) \subset 
\chainrec({\eps}) \subset \chainrec^{e} ({\sqrt{2n}\eps})$, 
the chain recurrent set $\chainrec = 
\cap_{\eps >0} \chainrec ({\eps})$ is 
the same for any metric uniformly equivalent to
euclidean.  Thus throughout we use $\chainrec ({\eps})$ as defined by 
our norm.  

\begin{rem}
When we say \textit{box}, we mean a ball around a point in this norm.  Note a box is also a vector of intervals, so boxes are neighborhoods which are easily manipulated with interval arithmetic.
%
\end{rem}

\subsection{The \boxchcn}
\label{sec:construction}
As suggested by the statement of Theorem~\ref{thm:babyRinBoxCov}, the \boxchcn is an inductive process.  We use the idea that any $\chreceps$ consists of precisely the $\eps$-pseudoperiodic orbits.  Below is an outline of the basic construction.

We start with a polynomial diffeomorphism $f$ of $\Ct$ (of dynamical degree $d(f)>1$), and let $F$ be any {\em interval extension} of $f$, i.e., $F$ is a function on interval vectors such that for any box $B$, the image $F(B)$ is a box containing $f(B)$ (see Appendix~\ref{sec:IA} for background on interval arithmetic).

\begin{enumerate}
\item[(0)] Given $f$ and an interval extension $F$, choose
 a small constant $\delta'_0 >0$, and a closed box $\boxcov_0$  in $\Ct$ which contains $\chrecdeltapo$ (of $f$).  Let $\eps_0$ be the side length of the box $\boxcov_0$.     Choose $\delta_0$ such that $\eps_0 \gg \delta_0 > \delta'_0$.
\item[(n)] Let $n \geq 0$. Suppose $\boxcov_n$ is a closed region in $\Ct$, consisting of a collection of boxes $\mathcal{V}_n = \{ B^n_k \}$ of side length at most $\eps_n$, and such that $\boxcov_n \supset \chrecdeltapn$, for some $\delta'_n >0$.  Suppose $\delta_n > \delta'_n$ is given (if $n\geq 1$, then $\delta_n$ is given by step ($n-1$)-(ii).)
\begin{enumerate}
\item[(i)]  Equally subdivide the boxes in $\mathcal{V}_n$.  That is,
choose $m > 1$, set $\eps_{n+1} = \eps_n / m$, and place a grid of $m^4$ boxes inside each box of $\mathcal{V}_n$.  This defines a new collection, $\mathcal{W}_{n+1}$, in which each box has side length $\eps_{n+1}$.
\item[(ii)] Build a graph approximating the map $f$ on $\mathcal{W}_{n+1}$.  Specifically, choose some $\delta_{n+1}$ such that $\delta_{n+1} < \delta_n/2$ and $\delta_{n+1} \ll \eps_{n+1}$ (for example, for every $n \geq 0$, set $\delta_{n} = \eps_{n} /1000$).
Then compute a directed graph $\Upsilon_{n+1}$ whose vertices are the boxes in $\mathcal{W}_{n+1}$, and
such that there is guaranteed to be an edge from box $B^{n+1}_k$ to box $B^{n+1}_j$ if $F(B^{n+1}_k)$ intersects a $\delta_{n+1}$-neighborhood of $B^{n+1}_j$.  
\item[(iii)] Find the subgraph of $\Upsilon_{n+1}$ consisting precisely of the vertices and edges which lie in cycles.  Call this subgraph $\Gamma_{n+1}$.   
Let $\mathcal{V}_{n+1} = \{ B^{n+1}_k \}$ be the vertices of $\Gamma_{n+1}$, and define
$\boxcov_{n+1}$ as the union of the boxes in $\mathcal{V}_{n+1}$.  
\end{enumerate}
\end{enumerate}

 %

\begin{rem} \label{rem:construction}
The \boxchcn immediately implies that statements (1) through (5) of Theorem~\ref{thm:babyRinBoxCov} are satisfied.
\end{rem}

\begin{rem}
The only difference between this basic \boxchcn and the procedure used by Osipenko (and others as discussed in the introduction) is the presence of the constants $\delta'_{n}$ and $\delta_n$.  In order to approximating $\chrec$, these constants are uneeded, and can all be taken to be zero.  However, in order to apply this construction to the problem of proving hyperbolicity, as we do in \cite{SLHtwo, SLHthree}, positive $\delta$'s are essential.
\end{rem}

The first step in verifying the usefulness of this construction is to show that if $\boxcov_{n-1} \supset \chrecdeltapnmo$, then there exists a $\delta'_{n}$ such that $\boxcov_{n} \supset \chrecdeltapn$.   This fact for polynomial diffeomorphisms of $\Ct$ follows from Theorem~\ref{thm:RinBoxCov} and Theorem~\ref{thm:babyRinBoxCov}, proved at the end of this section.


\subsection{Trapping Regions for $\chainrec$}
\label{sec:trap}

In order to apply the \boxchcn to polynomial diffeomorphisms of $\Ct$, we must first provide the base case, i.e., step (0) above.  For this we prove Proposition~\ref{prop:RinV}, in which we
calculate an explicit trapping region $\boxcov_0$ for the $\delta$-chain recurrent set of polynomial diffeomophisms of $\Ct$.  In particular, given a map $f$ and $\delta>0$, we give explicit $R'$ such that $\chrecdelta \subset \boxcov_0 = \{ \abs{x} \leq R', \abs{y} \leq R' \}$.

First we quantify how, for polynomial diffeomorphisms of $\Ct$,
infinity in the $x$ direction is attracting, while
infinity in the $y$ direction is repelling for $f$, and vice-versa for
$f^{-1}$. A version of the following lemma is given in \cite{FM}  (also see \cite{S1}).

\begin{lem} \label{lem:goodR}
Let $f$ be a polynomial diffeomorphism of $\Ct$, with $d(f)>1$.
Let $\delta >0$.
Then there is an $R>1$ and an $R' >R$, such that if $\sonorm{x}\geq R'$ and $\sonorm{x} \geq \sonorm{y}$, then 
$\snorm{f(x,y)} \geq \sonorm{x}+2\delta$.  

If $f_{a,c}$ is a \Henon mapping,  then
$R=\frac{1}{2}(1+\sonorm{a}+\sqrt{(1+\sonorm{a})^2 +4\sonorm{c}})$\
 and $2\delta = (R')^2 -(1+\abs{a})R' -\abs{c}$.
\end{lem}

\begin{proof}
Assume $f$ is a generalized \Henon mapping, $f(x,y) = (z,x)=(p(x)-ay,x)$,
with $d=$deg$(p)>1$.  So $p(x) = x^d + c_{d-1} x^{d-1} + \cdots + c_0$. 
Let
$q(r) = r^d - \sonorm{c_{d-1}} r^{d-1} - \cdots - \sonorm{c_0} -
(1+\sonorm{a})r$. Then there is an $R>0$ such that $q$ is monotone
increasing on~$[R,\infty)$, with $q(R)=0$. Note if $r>0$, then $q(r) \leq
r^d - r$, thus $R \geq 1$.

Since $q$ is a polynomial and is monotone increasing on $[R, \infty)$, with $q(R)=0$, we see $q : (R, \infty) \to (0, \infty)$ is invertible.  Thus given $\delta$, define $R' \in (R, \infty)$ by $q(R') = 2\delta$.

Let $(x,y)\in\CC$
satisfy $\sonorm{x}\geq R'$ and $\sonorm{x} \geq \sonorm{y}$.  Then 
\[
\sonorm{z} \geq \sonorm{p(x)}-\sonorm{a}\sonorm{y}
     \geq \sonorm{p(x)}-\sonorm{a}\sonorm{x}
    \geq q(\sonorm{x})+\sonorm{x}
    \geq 2\delta + \sonorm{x}.
\]
If $f = f_m \circ \cdots \circ f_1$, then each composition moves $\abs{z}$ farther away from $\abs{x}$ additively.  Thus, 
let $R'_k$ satisfy $q_k(R'_k) = 2\delta / m$ and take $R' = \max \{ R_k \}$, for $1 \leq k \leq m$. 

Note if $f_{a,c}$ is a quadratic \Henon mapping, $f_{a,c}(x,y)=(x^2+c-ay,x)$, then
$q(r) = r^2 - \abs{c} + (1+\abs{a})r$, hence we easily achieve the claimed bound. 
\end{proof}

\begin{defn} \label{defn:V}
Let $\delta >0$.  Following \cite{BS1}, we define the ``trapping regions'', for $R' \ $
given by Lemma~\ref{lem:goodR}, by:
\begin{eqnarray*}
\boxcov_0   & = & \{\sonorm{x} \leq R' \text{ and } \sonorm{y} \leq R'\}; \\ 
\boxcov_0^- & = & \{\sonorm{x} > R' \text{ and } \sonorm{x} > \sonorm{y}\}; \\
\boxcov_0^+ & = & \{\sonorm{y} > R' \text{ and } \sonorm{y} > \sonorm{x}\}. 
\end{eqnarray*}
Note the regions $\boxcov_0 = \boxcov_0 (\delta)$ depend on $\delta$, but in each instance the $\delta$ will be clear from context, so we will rarely use the notation $\boxcov_0(\delta)$, rather we simply use $\boxcov_0$.
\end{defn}

The sets $V, V^{\pm}$ introduced in \cite{BS1} are equal to $\boxcov_0, \boxcov_0^{\pm}$ for $\delta=0$, and satisfy 
$
K^+ \subset V \cup V^+, \ \ 
K^- \subset V \cup V^-, \ \text{ and } 
K\subset V$ (\cite{BS1}). 
Choosing $R'$ larger than $R$ allows us to preserve these relationships and trap $\delta$-pseudo-orbits
as well.

\begin{prop} \label{prop:RinV}
Let $f$ be a polynomial diffeomorphism of $\Ct$, with $d(f)>1$.
Let $\delta>0$ and let $\boxcov_0$ be as in Definition~\ref{defn:V}.
Then
$
\chrecdelta \subset \boxcov_0.
$
\end{prop}

\begin{proof}
Given Lemma~\ref{lem:goodR}, we have
$f(\boxcov_0^-) \subset \boxcov_0^-$ and $f(\boxcov_0^-) \cap 
\snbd{\boxcov_0}{2\delta} = \emptyset$.  
Thus if $p \in \boxcov_0^-$, then $p$ is not in
$\chrecdelta$, since 
the images
$f(x_k)$ move by at least $2\delta$ in the $x$ direction, and so
the $x_{k+1}$ coming back in by only $\delta$ makes it impossible
for $x_n=p$.  

Similarly,  for $p \in \boxcov_0^+$ look at the chain backwards to contradict $\delta$-chain recurrence.
\end{proof}

To get an idea of the size of $R$, note that  $R=2$ for
$c=2$.  Since the Mandelbrot set is contained in $\{ x : \abs{x}
\leq 2\}$, for the parameters we tend to study we have
$1\leq R\leq 2$.  For a \Henon mapping, the values of $R$ are also close to
this range.

\subsection{Defining the graphs $\Upsilon$ and $\Gamma$}
\label{sec:graphs}

In steps (ii) and (iii) of the \boxch construction, we compute graphs $\Upsilon$ and $\Gamma$ representing the action of $f$ (or $F$) on our collection of boxes.  The following terminology will ease our discussion of these graphs.

About notation:  if $\Gamma$ is a graph, then $\mathcal{V}(\Gamma)$ denotes the vertex set of $\Gamma$, and $\mathcal{E}(\Gamma)$ denotes its edge set.
Also we often discuss one subset of a collection of boxes $\mathcal{V}  = \{ B_k \}_{k=1}^N$ at a time, and since the ordering is unimportant we avoid double subscripts and simply use $\{ B_0, B_1, \ldots, B_{\ell}\}$.  

\begin{defn} \label{defn:boxmodel}
Let  $\Upsilon = \Upsilon_n$ be the directed graph built in step (n)-(ii) of a \boxch construction.
Then we say $\Upsilon$ is 
an \textit{$(\eps, \delta)$-\boxchmod of $F$}.

In addition, when we say $\Upsilon$ is an  $(\eps, \delta)$-\boxchmod of $f$, we mean the ``theoretically ideal'' model, using the ideal interval extension of $f$, \textit{i.e.}, require $F(B_k)=$ Hull$(f(B_k))$, for all boxes $B_k$.  
\end{defn}

Note for any interval extension $F$ of $f$,
 we know $F(B) \supset \text{Hull}(f(B))$ for any box $B$. Thus any result which is true for all interval extensions $F$ of $f$ is also true for $f$.  Thus our default is to discuss box chain models of $F$, and only use box chain models of $f$ when trying to be precise about theoretical estimates.  


Suppose $\Upsilon$ is an $(\eps, \delta)$-\boxchmod of $F$.  Note by Definition~\ref{defn:boxterms}, the subgraph $\Gamma$ consisting of the cycles of $\Upsilon$ is called an $(\eps, \delta)$-\boxchrecmod of $F$.  

\begin{rem}
Both a \boxchmod and a \boxchrecmod of $f$ satisfy the definition of a {\em symbolic image of $f$}, given by Osipenko in \cite{Osiold}.
\end{rem}

When the context is clear, or the distinction is unimportant, we simply refer to a box chain recurrent model as a box chain model.  Thus we use the symbol $\Upsilon$ for any graph which is either a box chain model or a box chain recurrent model, and reserve $\Gamma$ for box chain recurrent models.

The following standard concept in graph theory will help us to analyze $\Gamma$.

\begin{defn} \label{defn:scc}
Let $\Upsilon$ be a directed graph.  If there is a path from vertex $v$ to 
vertex $u$, then we say $u$ is reachable from $v$.  
A {\em strongly connected component} (SCC)
is an equivalence class under the ``are mutually reachable'' equivalence 
relation.  If an SCC 
consists of only one vertex, it must have an edge to itself.  
\end{defn}

Note the following easy relationship.

\begin{lem}
Let $\Upsilon$ be any directed graph. Let $\Gamma$ be the subgraph of $\Upsilon$ consisting of the vertices and edges which lie in cycles.
Then $\Gamma = ( \Gamma^1 \sqcup \Gamma^2 \sqcup \ldots \sqcup \Gamma^l )$ is precisely the subgraph of $\Upsilon$ consisting of the union of the SCC's of $\Upsilon$. 
\end{lem}

Recall from Definition~\ref{defn:boxterms} that the {\em \boxchtrcomps} of $\Gamma$ are defined as the edge-connected components of $\Gamma$, hence these are precisely the SCC's of  $\Upsilon$.
The relationship between a \boxchtrcomp and an $\eps$-chain transitive component is made precise in Corollary~\ref{cor:sccchaincomp}, and justifies our terminology.
Thus the decomposition of a graph $\Gamma$ into its SCC's is analogous to the partitioning of the chain recurrent set into it's invariant pieces, the chain transitive components.  

 \cite{TCR} gives a standard procedure for decomposing a graph into its SCC's.
 

\subsection{Estimates on  $\Upsilon$ and $\Gamma$ 
}
\label{sec:modelbounds}

Now we can calculate estimates which quantify our models.   
%
%
%
%
%
We begin with a lemma on the size of the image of the boxes.

\begin{lem} \label{lem:imageboxsize}
Let $\Upsilon_n$ be an $(\eps_n, \delta_n)$-\boxchmod of $f$. 
Then there exists $r_n>0$ (depending on $\eps_n$, $f$, and $\boxcov_0$) such that for $\eps_n' =(r_n + 1)\eps_n+\delta_n$,
and for any $B_k \in \mathcal{V} ({\Upsilon_n})$, we have:
\begin{enumerate}
\item the side length of the box $\text{Hull}(f(B_k))$ is less than or equal to
  $\eps_n r_n$, and
\item if $(k,j) \in \mathcal{E}({\Upsilon_n})$, then
for any $x_k \in B_k$ and any $x_j \in B_j$, 
$\snorm{f(x_k) - x_j} < \eps'_n$. 
\end{enumerate}
For $f_{a,c}$ a \Henon mapping, we may take
$r_n=\eps_n+ (2R'+\sonorm{a})$, where
$R'$ is as in Proposition~\ref{prop:RinV}.
\end{lem}

\begin{proof}
The second item follows immediately from the first item, and the fact
that if
$B_j \cap \snbd{\text{Hull}(f(B_k))}{\delta} \neq \emptyset$, then there must be an edge from $B_k$ to $B_j$.
We prove the first item using the linearization of $f$ to
approximate it.

We assume $f$ is a generalized \Henon mapping, so $f(x,y) = 
(p(x)-ay,x)$, $p$ monic of degree $d>1$.
Let $B \in \mathcal{V}({\Upsilon})$ and $(z,w) \in B$.  
The linearization of $f$ at $(z,w)$ is
\[
L_z f \twovec{x}{y}  =  f \twovec{z}{w} + D_zf \twovec{x-z}{y-w} 
=   \twovec{p(z)-p'(z)(x-z)-ay}{x}.
\]
Note $L_z f(z,w) = f(z,w)$, for any $w$. 
Next, we observe
\[
\snorm{f \twovec{x}{y} - L_z f \twovec{x}{y}} 
 =  \snorm{\twovec{p(x)-p(z)-p'(z)(x-z)}{0}} 
 \leq  \sum_{k=2}^d \frac{\abs{p^{(k)}(z)}}{k!}  \abs{x-z}^k.
\]
If $(x,y)$ is also in~$B$, then $\sonorm{x-z} \leq \eps_n$.  Since $p$ is
a polynomial, and $z\in \boxcov_0=\snbd{0}{R'}$, then for any $k \geq 0$, there 
exist $T_k \geq 0$ such that $\abs{p^{(k)}(z)}/k! \leq T_k$. Hence,
\[
\sum_{k=2}^d \frac{\abs{p^{(k)}(z)}}{k!}  \abs{x-z}^k
\leq \sum_{k=2}^d T_k \eps_n^k.
\]
Now we need to bound
\[
\snorm{L_z f\twovec{x}{y} - L_z f \twovec{z}{w}}
 =  \max\(\sonorm{x-z}, \sonorm{p'(z)(x-z)-a(y-w)}\).
\]
If $(x,y)$ is in~$B$, then we also know $\sonorm{x-z} \leq \eps_n$, and
\[
\sonorm{p'(z)(x-z)-a(y-w)}  \leq  T_1 \eps_n + \sonorm{a} \eps_n.
\]

Finally, we put the above pieces together to compute, for any $(x,y),
(z,w)$ in~$B$,
\[
\snorm{f \twovec{x}{y} - f\twovec{z}{w}}
 \leq 
 \sum_{k=2}^d T_k \eps_n^k + \eps_n \max \(1, T_1 +\sonorm{a} \).
\]

Hence, if $r_n$ is set to be the maximum above divided by $\eps_n$, then we have that the diameter of
the set $f(B_k)$ in the metric $\snorm{\cdot}$ is less than or equal to
$r_n\eps_n$, hence the side length of Hull($f(B_k)$) is less than or equal to $r_n \eps_n$.

For $f_k$ a generalized \Henon mapping, we have $r_n(f_k)$ as defined above.  Now, if $f= f_m \circ \cdots \circ f_1$, then set $r_n=r_n(f_1) \cdots r_n(f_m)$.
This suffices, for if we let $s(f_k,\eps)=r_n(f_k) \eps$ be the bound on the diameter of a box of size $\eps$ under $f_k$, then we get $s_1 = s(f_1,\eps_n) =  \eps_n r_n(f_1),$ then  $s_2 = s (f_2 \circ f_1, \eps_n)  = s(f_2, s_1) = s_1 r_n(f_2) = \eps r_n(f_1) r_n(f_2),$ etc., and 
$s_m = s(f_m \circ \cdots \circ f_1) =  s(f_m, s_{m-1}) = \eps r_n(f_1) \cdots r_n(f_m).$

In the case of $f_{a,c}$ a \Henon mapping, since $p(z)=z^2+c$ we can compute that $T_1=2 R'$ and $T_2=1$, and recall from Proposition~\ref{prop:RinV}, we have 
$R'>1$. Hence $r_n=\eps_n + \max \( 2 R' + \sonorm{a},1)=\eps_n + 
(2R'+\sonorm{a}\)$.
\end{proof}

This lemma give us the following component of Theorem~\ref{thm:babyRinBoxCov}.

\begin{cor} \label{cor:imageboxsize}
Item (6)-(a) of Theorem~\ref{thm:babyRinBoxCov} is satisifed. 
\end{cor}

\begin{proof}
First, we show there exists a $C>0$ such that for any $n > 1$, $r_n \leq C-1$, hence $\eps'_n \leq \delta_n + C\eps_n$.  

Suppose for all $n\geq 1$ that $\eps_n < 1$. Then we know $\eps_n^k < \eps_n$ for any $k > 1$. Define $C= 1+ \sum_{k=2}^d T_k + \max \(1, T_1 +\sonorm{a} \)$ for a generalized \Henon mapping, or in the manner analogous to the previous proof for a composition (\textit{i.e.}, for $f = f_m \circ \cdots \circ f_1$, take $C = (C_m -1) \cdots (C_1 -1) + 1$).  Then we see $r_n \leq C-1$.  

On the other hand, since $\eps_n$ is a decreasing sequence (see Remark~\ref{rem:construction}), at most  $\eps_n >1$ for $0 \leq n \leq N$.
Thus if we set $C = 1+ \max \{ r_1, \ldots r_N, \sum_{k=2}^d T_k + \max \(1, T_1 +\sonorm{a} \) \}$, then $r_n \leq C-1$ for all $n\geq 1$. 

Finally, the above lemma defines $\eps'_n$ so that $\eps'_n > \eps_n$. Hence we have $\eps_n < \eps'_n \leq \delta_n + C\eps_n$.  By construction (see Remark~\ref{rem:construction}) we know both $\eps_n$ and $\delta_n$ decrease to zero as $n \to \infty$, hence $\eps'_n$ must as well. 
\end{proof}

The next lemma will allow us to find $\delta'$ for which a \boxchrecmod traps all $\delta'$ pseudo-periodic orbits, \textit{i.e.}, $\chrecdeltap$.

\begin{lem} \label{lem:edgeerror}
Let $\Upsilon_n$ be an $(\eps_n, \delta_n)$-\boxchmod of $F$.
There exists an $\eta = \eta(n) \in (0,\delta_n)$, such that
if $x \in \snbd{B_k}{\eta}$ and 
$f(x) \in \snbd{B_j}{\eta}$, then there is an 
edge from $B_k$ to $B_j$ in~$\Upsilon_n$, \textit{i.e.}, $(k,j)
\in \mathcal{E}({\Upsilon_n})$.

For $f_{a,c}$ a \Henon mapping, we can take
\[
\eta(n) = \frac{1}{2} \( -(2 R'+\sonorm{a}+1) + 
 		\sqrt{(2 R'+\sonorm{a}+1)^2+4\delta_n} \),
\]
where $R'$ is as in Proposition~\ref{prop:RinV}.
\end{lem}

\begin{proof}
Again we give the proof for $f(x,y)=(p(x)-ay,x)$,
$p$ monic of degree $d>1$.
Let $(x,y) \in 
\snbd{B_k}{\eta}$ and $f(x,y)  \in 
\snbd{B_j}{\eta}$.  Let $(z,w)$ be a point in $B_k$ which 
realizes this minimum distance (since $B_k$ is closed), \textit{i.e.},
$
\snorm{(x,y) - (z,w)} < \eta.
$
Then in order to guarantee that $(k,j) 
\in \mathcal{E}({\Upsilon})$, we 
just need
$
\snorm{f(x,y) - f(z,w)} < \delta_n - \eta.
$
But examining the proof of Lemma~\ref{lem:imageboxsize}, since $(z,w) \in 
B_k$, we see that we have
\[
\snorm{f\twovec{x}{y} - f\twovec{z}{w}} \leq 
\sum_{k=2}^d T_k \eta^k + \eta \max \(1, T_1 +\sonorm{a} \).
\]
Thus, we just need $\eta$ to satisfy
\[
\sum_{k=2}^d T_k \eta^k + \eta \(\max \(1, T_1 +\sonorm{a}\)+1\)  - \delta_n 
<0.
\]
Let $q_n(t) = \sum_{k=2}^d T_k t^k + t \max \(2, T_1 +\sonorm{a}+1 
\) - \delta_n$.
Then we set $\eta = \eta(n)$ to be the smallest positive root of $q_n$.   
If $f=f_m \circ \cdots \circ f_1$, we may take $\eta = \min\{ \eta (f_1), \ldots,
\eta (f_m)\}$. For $f_{a,c}$ a \Henon mapping, we get $q_n(t) = t^2 + 
t(2R'+\sonorm{a}+1)-\delta_n$, which leads to the claimed bound.
\end{proof}

\begin{cor} \label{cor:deltap}
Recall $\boxcov_0$  denotes the initial bounding box, with $\chrecdeltapo \subset \boxcov_0$ for some $\delta'_0 > 0$.  For each $n\geq 1$, set $\delta'_n = \min( \eta(n), \delta'_0)$, with $\eta(n)$ as in Lemma~\ref{lem:edgeerror}.

Then item (6)-(b) of Theorem~\ref{thm:babyRinBoxCov} is satisfied.  
\end{cor}

\begin{proof}
Lemma~\ref{lem:edgeerror} states that $\eta(n) < \delta_n$, hence $\delta'_n < \delta_n$.  Next, note from the proof of Lemma~\ref{lem:edgeerror} that
$\eta$ is a monotone increasing function of $\delta_n >0$.  By Remark~\ref{rem:construction} the sequence $\delta_n$ decreases as $n$ increases. Hence as $n$ increases,  we get $\eta(n)$ decreases.  Thus $\delta'_n$ is at most equal to $\delta'_0$ for a finite number of terms $N$, and then for $n > N$,  $\delta'_n$ decreases.
\end{proof}

\subsection{Quantifying the accuracy of an $(\eps, \delta)$-\boxchrecset $\boxcov$}
\label{sec:resboxcov}

In this section, we combine the lemmas of the previous section to achieve Theorem~\ref{thm:RinBoxCov}, which immediately implies Theorem~\ref{thm:babyRinBoxCov} along with several relevant corollaries.

\begin{thm} \label{thm:RinBoxCov}
Suppose 
$\Upsilon_n$ is an $(\eps_n, \delta_n)$-\boxchmod of $f$, and
$\Gamma_n$ is the $(\eps_n, \delta_n)$-\boxchrecmod of $f$ consisting of the SCC's of $\Upsilon_n$.  
Let $\eps'_n$ be as in Lemma~\ref{lem:imageboxsize} and
let $\delta'_n$ be as in Corollary~\ref{cor:deltap}.  

Let $\boxcov (\Gamma_n) = \boxcov_n$ denote the region in $\Ct$ covered by the box vertices of $\Gamma_n$. Then item (6)-(c) of Theorem~\ref{thm:babyRinBoxCov} is satisfied, \textit{i.e.},
\[
\chrecdeltapn \subset \boxcov (\Gamma_n) \subset \chrecepspn.
\]
\end{thm}

\begin{proof}
We use the fact that precisely the vertices which lie in cycles in $\Upsilon_n$ are those which are in some strongly connected component of $\Upsilon_n$.  

We handle the two inclusions separately. First consider the inclusion $\chrecdeltapn \subset \boxcov (\Gamma_n)$. We first need to establish $\chrecdeltapn \subset \boxcov(\Upsilon_n)$,
which we  prove by induction. 
For $n=0$, we can choose any $\delta'_0 >0$, and produce the box $\boxcov_0$ such that $\chrecdeltapo \subset \boxcov_0$ by Proposition~\ref{prop:RinV}.  
Consider $\Upsilon_0 = \Gamma_0$ to be the graph with the single vertex $\boxcov_0$ with an edge to itself.
Now suppose $n \geq 1$, and assume 
$\chrecdeltapnmo \subset \boxcov (\Gamma_{n-1}) $. 
Observe that since the boxes of $\Upsilon_n$ are obtained from subdividing the boxes of $\Gamma_{n-1}$, we have $\boxcov(\Upsilon_n) = \boxcov(\Gamma_{n-1})$, hence   
$\chrecdeltapnmo \subset \boxcov(\Upsilon_n)$.  
Also, 
since $\delta'_n \leq \delta'_{n-1}$ (Corollary~\ref{cor:deltap}) we have $\chrecdeltapn \subset \chrecdeltapnmo$. Thus
$\chrecdeltapn  \subset  \boxcov({\Upsilon_n})$. 

Suppose $p \in \chrecdeltapn$.  Then
there exist $x_1=p, x_2, \ldots, x_{m-1}, x_m=p$ such that $\snorm{f(x_k)-x_{k+1}}<
\delta '_n$ for $1\leq k < m$.  Note that $x_k \in \chrecdeltapn$,
for $1\leq k < m$. Hence each $x_k \in \boxcov ({\Upsilon_n})$  as well.  Then there are
boxes $B_k\in \Upsilon_n$ such that $x_k \in B_k$ for $1\leq k < m$.  Since
$\snorm{f(x_k)-x_{k+1}}< \delta'_n$, we have $f(x_k) \in
\snbd{B_{k+1}}{\delta'_n}$.  Since $\delta'_n \leq \delta_n$, there is an edge
in $\Upsilon_n$ from $B_k$ to $B_{k+1}$.

Hence, $p$ is in a box $B_1$ which lies in a cycle of $\Upsilon_n$, $B_1 \to
B_2 \to \cdots \to B_{m-1} \to B_1$.  Thus $B_1 \in \mathcal{V}({\Gamma_n})$,
hence $p\in \boxcov_n$.

For the second inclusion, $\boxcov (\Gamma_n) \subset \chrecepspn$,
suppose $p\in \boxcov (\Gamma_n)$.  Thus $p$ lies in some box $B_1$
 which lies in
a cycle $B_1\to B_2\to \cdots \to B_{m-1} \to B_1$ in $\Upsilon_n$.  Recall
$\eps_n$ is the side length of the boxes in the graph $\Gamma_n$.  
  Let
$x_1=x_m=p$, and $x_k$ be any point in~$B_k$ for $2\leq k \leq m-1$. Then
by Lemma~\ref{lem:imageboxsize}, since $B_1\to B_2\to \cdots B_{m-1} \to
B_1$ is a cycle in $\Upsilon_n$, we have $\snorm{f(x_k)-x_{k+1}} < \eps'_n$
for $1\leq k\leq m$.  Hence, $p$ is $\eps'_n$-chain recurrent.
\end{proof}


Note this theorem completes the proof of Theorem~\ref{thm:babyRinBoxCov}.

A first modification of this theorem is that if we consider $\delta'=0$ in the hypothesis, we can conclude $\chainrec \subset \boxcov_n \subset \chrecepsp$.  We need this observation in Section~\ref{sec:codetricksone} to justify the process of eliminating ``$\boxcov_0$-escaping'' boxes.

Next, extending from the theoretical $f$ to the practical $F$, we immediately obtain:

\begin{cor} \label{cor:RinBoxCov}
Suppose $\Upsilon_n$ is an $(\eps_n, \delta_n)$-\boxchmod of $F$, and
$\Gamma_n,  \eps'_n, \delta'_n,$ and  $\boxcov_n$ are as in Theorem~\ref{thm:RinBoxCov}.  

Let
\[
s_n({F}) = \sup_{B_k \in \mathcal{V}(\Gamma_n)}  \text{sidelength} [ F(B_k) - \text{Hull}(f(B_k)) ].
\]

Then
\[
\chrecdeltapn \subset \boxcov_n \subset \chrec(\eps'_n+s_n(F)).
\]
\end{cor}

We cannot compute $s_n(F)$ exactly, but simply strive to design an implementation to minimize it.  We shall not discuss that matter in detail in this paper, but refer the reader to the interval arithmetic resources listed in Section~\ref{sec:IA}.  Away from the boundaries of machine precision, it is reasonable to assume $s_n(F)$ decreases as $\eps_n$ decreases, and that $s_n(F)$ is small compared to $\eps'_n$.  In fact, for a naive interval extension of $f$, $s_n(F) = O(\eps_n)$.

We can also immediately apply the theorem to the SCC's.

\begin{cor} \label{cor:sccchaincomp}
Assume the hypothesis of Theorem~\ref{thm:RinBoxCov}.
\begin{enumerate}
\item Let $\chrecpdeltap$ be any $\delta '$-chain
transitive component.  Then there is a \boxchtr component, $\Gamma'$,
such that $\chrecpdeltap \subset \boxcov ({\Gamma'})$. 
\item Let ${\Gamma'}$ be any \boxchtrcomp of $\Gamma$. Then 
there is some $\eps'$-chain transitive component, $\chrecpepsp$, 
such that $\boxcov({\Gamma'}) \subset \chrecpepsp$.
\end{enumerate}
\end{cor}

Again, in Corollary~\ref{cor:sccchaincomp} with the weaker hypothesis $\delta'=0$, we can form a conclusion analogous to (1) by setting $\delta'=0$, \textit{i.e.}, $\chrecdeltap = \chainrec$.

Note that more than one $\delta'$-chain transitive component
can be contained in a single \boxchtr component.  Also, an
$\eps'$-chain component may not actually contain any chain recurrent
points, hence a strongly connected component may not contain any chain
recurrent points. We explore this further in Section~\ref{sec:sinkJsep}.

For polynomial diffeomorphisms of $\Ct$ (with $d(f)>1)$, the Julia set $J$ is contained in a 
single chain transitive component $\chainrec'$ of $\chainrec$ (see Theorem~\ref{thm:JinR}).
Hence Corollary~\ref{cor:sccchaincomp} implies:

\begin{cor} \label{cor:JinB}
Assume the hypothesis of Theorem~\ref{thm:RinBoxCov}.
Then there is a single \boxchtrcomp 
${\Gamma'}$ such that  $J \subset  \boxcov({\Gamma'})$.
\end{cor}

Thus, one of the \boxchtr components, $\Gamma'$, contains $J$.  The others either 
contain attracting (or repelling) periodic orbits, or do not intersect 
$\chainrec$.  We can easily identify the component containing $J$, for it has by far the most vertices.

\begin{rem} \label{rem:dimone}
Note the \boxchcn and all of the estimates and results of this section can be reformulated to apply to polynomial maps of $\CC$ of degree $d>1$, by simply dropping the $y$-component, and setting $a=0$ in all of the calculations.  We use these results to study polynomial maps of $\CC$ in \cite{SLHtwo}.
\end{rem}

See Section~\ref{sec:examples} for examples of \boxchrecmods generated with the program Hypatia.


\section{Improving efficiency using Dynamics}
\label{sec:codetricksone}

One result of our interest in the complex case is that we are forced to deal with increased computational complexity difficulties, since 
working in $\Ct$ is the same for a computer as working in $\RR^4$.  In our experience, the main computational limitation is memory usage 
(even using a computer with 4 GB of RAM).  Thus we keep memory efficiency in mind when tailoring the \boxchcn  to polynomial diffeomorphisms of $\Ct$.
In this section, we discuss how to insert into the \boxchcn  two improvements, designed to increase efficiency, but without completely negating Theorem~\ref{thm:RinBoxCov}.  In designing these algorithms, we take advantage of the dynamics of the map.

\subsection{Step (i$'$): Selective Subdivision}

Note first that a single step (rather than inductive) construction consisting of simply subdividing the initial $\boxcov_0$ into a very large grid of very small boxes is conceptually much simpler, but very computationally inefficient.
%
%
In addition to the obvious work that subdivision saves, the tree
structure it creates is useful for quickly computing things
such as which boxes intersect a given set (like the image of another box). 
This allows us to store each graph as an array of vertices, with edges as
adjacency lists, with no need to arrange the vertices in any particular
order in the array.
Further, the tree structure created by subdivision lends itself easily to be improved in the following manner.

In trying to approximate a dynamically defined set, it is natural to improve upon the basic subdivision procedure by replacing step (i) of the \boxchcn  with the more sophisticated \textit{selective subdivision} procedure. 
The philosophy here is that we can concentrate our resources on the regions which are most ``troublesome'' by allowing boxes of different sizes.
We would ideally define grid boxes so that the dynamical behavior in each box is varying by at
most a small amount. Thus, we would like to refine down to a certain
reasonable box size, then somehow select only a small fraction of the
boxes to be subdivided further, leaving the rest unchanged. 

\begin{enumerate}
\item[(i$'$)]  Sort the boxes in $\mathcal{V}_n$ into two sets: $\mathcal{V}^s_n$ is boxes to be subdivided, and $\mathcal{V}^u_n$ is boxes to remain unchanged.  Equally subdivide the boxes in $\mathcal{V}^s_n$: 
choose $m > 1$  and place a grid of $m^4$ boxes inside each box of $\mathcal{V}^s_n$, to obtain the collection $\mathcal{W}^s_{n+1}$.  Set $\mathcal{W}_{n+1} = \mathcal{W}^s_{n+1} \cup \mathcal{V}^u_n$. Then each box in $\mathcal{W}_{n+1}$ has side length at most $\eps_{n+1} = \eps_n$, and at least $\eps_n/m$.
\end{enumerate}

Note in this case the sequence of maximum box side lengths $\{ \eps_n \}$ is nonincreasing, rather than decreasing.  In order to acheive the ideal accuracy in the limit, as in Theorem~\ref{thm:RinBoxCov}, one would have to subdivide all of the boxes once every few steps.

To implement this procedure, we must somehow choose which boxes to subdivide.
A first goal is to obtain a \boxchrecmod which
separates the Julia set from the attracting and repelling periodic orbits.  With the goal of ``refining out'' the attracting behavior, we implemented in Hypatia an option to subdivide 
only boxes which seem to be in $K^{+}$. We call this \textit{\inKsel subdivision}.
For example, one way of detecting boxes in sink basins is to choose boxes such that all eigenvalues of a few iterates of the
derivative matrix $Df$ are small.  
 In some cases, \inKselsub does help to separate out the
sink dynamics more quickly (see Example~\ref{exmp:altpertwo}).
However, there are interesting examples for which this still does not allow us to separate the sink from $J$ (see Example~\ref{exmp:per31}).
Thus, it seems that a selective subdivision procedure could be an extremely useful step, however it is unclear what are the optimal selection criteria for a given map.

\subsection{Step (i.5): Eliminating $\boxcov_0$-escaping boxes}
After performing step (i) or (i$'$), \textit{i.e.}, subdividing the desired level $n$ boxes $\mathcal{V}_n$ to obtain a new collection of level ${n+1}$ boxes, $\mathcal{W}_{n+1}$, but before performing step (ii) (\textit{i.e.}, computing a graph representing the action of $f$ (or $F$) on the collection of boxes), we institute an additional efficiency improving check:  \textit{eliminating $\boxcov_0$-escaping boxes}.
This is a very useful technique which removes much of the work of computing
the graph in step (ii).  The idea is to eliminate boxes whose images eventually lie outside of~$\boxcov_0$.

Note by examining the proof of Proposition~\ref{prop:RinV} that if  $f^m(B_k) \cap \boxcov_0 = \emptyset$ for some $m >1$, then $B_k$ contains no points of 
$\chainrec$.  
Thus we can delete vertex $k$ from the graph
before computing the edge list.  For the \Henon family, we can also take
advantage of invertibility.  
That is, it follows from
Lemma~\ref{lem:goodR} that if $f^{-m}(B_k) \subset \CC \setminus \boxcov_0$, then
$B_k$ contains no points of 
$\chainrec$.  Thus we can check forward and
backward images of boxes, deleting any boxes with some image leaving $\boxcov_0$,
and then in step (ii) compute the edge list for the reduced graph.

\begin{enumerate}
\item[(i.5)]  Given a collection of boxes $\mathcal{W}_{n+1}$, eliminate the $\boxcov_0$-escaping boxes.  Remaining is a subcollection of boxes $\mathcal{W}^h_{n+1}$ (which contains $\chainrec$).  Replace $\mathcal{W}_{n+1}$ with $\mathcal{W}^h_{n+1}$.
\end{enumerate}

Each $\boxcov_n$ in the constructed sequence contains $\chainrec$, so we may call it an $(\eps, 0)$-approximation to $\chainrec$.
Such a sequence satisfies slightly weaker versions of
Theorem~\ref{thm:RinBoxCov} and Corollary~\ref{cor:sccchaincomp}, with
${\delta'}$ replaced by $0$, and $\chrecdeltap$ 
replaced by $\chainrec$.  Note however we still use the $\delta$ factor to build the edges in the graph, thus we still satisfy Lemma~\ref{lem:imageboxsize} and Lemma~\ref{lem:edgeerror} (this is important in \cite{SLHtwo} and \cite{SLHthree}).

In many instances, this check
eliminates half or even three-fourths or more of the new boxes.
Since the main computational limitation of {\Hypatia} is in memory usage, not
having to store edges corresponding to boxes which are eliminated with
this iteration check is a big gain.
Specific data for some examples is in Table~\ref{table:examples}, in Section~\ref{sec:examples}.


\section{Separating $J$ from a fixed sink}
\label{sec:sinkJsep}

In this section, we calculate a theoretical limitation of the \boxchcn in the case of 
a \Henon mapping, $f(x,y) =
(x^2+c-ay, x)$, with an attracting fixed point $p=(z,z)$.  
In particular, we establish Theorem~\ref{thm:babysinkJsep},  
quantifyiing a box size which guarantees that the \boxchcn produces 
a \boxchrecmod {\em separating} the fixed sink from $J$, as in Definition~\ref{defn:sep}, \textit{i.e.}, 
such that the fixed sink lies in a separate \boxchtrcomp from $J$.  This 
quantification is in terms of $a,c$, and the eigenvalues, $\lambda_1 
\neq \lambda_2$, of $D_pf$, and $\lambda = \max (\abs{\lambda_1}, \abs{\lambda_2})$.

First we find a euclidean disk contained in the
sink basin, in Proposition~\ref{prop:sinkbasin}.
Second, we quantify an annular region about the sink which contains only non 
$\eps$-chain recurrent points, in Proposition~\ref{prop:babynonepsrec}.
We use this along with the estimates of Section~\ref{sec:boxcover} to derive 
Theorem~\ref{thm:babysinkJsep}.
Finally, we apply our estimates to the 3-1-map.

\begin{exmp}[The 3-1-map] \label{exmp:per31}
Recall from Example~\ref{exmp:introper31} that the \Henon mapping $f_{a,c}$ with $(a,c) = (.3, -1.17)$  
 is an interesting example
because it appears to have two attracting periodic cycles, one of
period three and one of period one.  Two attracting cycles is not a phenomenon which can
occur for the quadratic polynomial $f_c(z) = z^2+c$.  Unfortunately, using our program \textit{Hypatia} implementing the \boxch construction, we failed to separate the sinks from $J$ before running out of the 
$4$ GB of RAM available on our computer.

In an attempt to find a good \boxchrec set, we first uniformly
subdivided all boxes to obtain a $(2^7)^4$ grid on $\boxcov_0 = \snbd{0}{2.01}$, with box side
length $2R/2^7 = 0.03$, where $R=2.01$. Then we used \inKselsub (Section~\ref{sec:codetricksone}), which subdivided about half
of the boxes.  At this point, the
smallest boxes had side length $2R/2^8=0.015$.  The \boxchrecmod  $\Gamma$
was composed of~$944{,}000$ boxes and $66{,}500{,}000$ edges.  This used
approximately $3.2$ GB of RAM, thus it seemed we could not subdivide
significantly farther.

We also tried uniformly subdividing to obtain a $(2^6)^4$ grid on 
$\boxcov_0$, then 
invoking \inKsel subdivision, twice, to get some boxes as small as above.  However, 
this did not significantly decrease the amount of memory used.

Figure~\ref{fig:per31} shows the unstable manifold slice of the 
\boxchrecset from the $(2^7)^4$ grid on $\boxcov_0$.
We were initially surprised that we could not achieve separation for this map.  
This motivated the estimates of this section. 
\end{exmp}

\subsection{A dynamically significant norm}

To quantify the dynamical notations of interest, we need to start with a
metric which respects the dynamics.

\begin{defn} \label{defn:vnorm}
Note the Jacobian of $f$ is
\[
D_{(x,y)} f = \begin{bmatrix}
2x & -a \\ 1 & 0
\end{bmatrix}
\] 
Suppose $\lambda_1 \neq \lambda_2$ are the eigenvalues of
$D_pf$, for the fixed sink $p=(z,z)$.  Let $\lambda = \max( \abs{\lambda_1},\abs{\lambda_2})$. Note
that since $p$ is a fixed sink, $\abs{\lambda} < 1$.  Let 
$\{\mathbf{v_1}, \mathbf{v_2}\}$ be the basis of eigenvectors, where we
choose $\mathbf{v_j} = (\lambda_j, 1)$.  Then $\mathbf{A}= [ \mathbf{v_1}
\ \mathbf{v_2} ]$ is the change of basis matrix, \textit{i.e.}, if $\{ 
\mathbf{e_1},
\mathbf{e_2} \}$ is the standard basis in $\Ct$, then
$\mathbf{A}\mathbf{e_j} = \mathbf{v_j}$.  Let $\enorm{\cdot}$ be the 
euclidean norm in $\Ct$.  Define the norm $\vnorm{\cdot}$ by
$
\vnorm{\mathbf{u}} := \enorm{\mathbf{A}^{-1}\mathbf{u}}.
$
\end{defn}


We show below that $f$ is 
contracting with respect to $\vnorm{\cdot}$ in an neighborhood of~$p$.  
%
%
First, we show this metric is uniformly equivalent to euclidean, and 
compute the constants of equivalence.

\begin{lem} \label{lem:vnormequiv}
For all $\mathbf{u} \in \Ct$,
$
C \vnorm{\mathbf{u}} \leq \enorm{\mathbf{u}} \leq D \vnorm{\mathbf{u}},
$
where $C,D$ are positive constants given by
\[
C = \frac{ \abs{\lambda_1 -\lambda_2}}
{\sqrt{2+\abs{\lambda_1}+\abs{\lambda_2}}},
\ \ \ \
D = \sqrt{2+\abs{a} + \lambda^2}.
\]
\end{lem}

\begin{proof}
Let $(x,y)$ be any vector in $\Ct$.
To calculate both C and D, we use the following observation: 
\begin{equation} \label{eqn:crossest}
0 \leq (\abs{x} - \abs{y})^2 \text{ \ implies \ }
2\abs{x} \abs{y} \leq \abs{x}^2 + \abs{y}^2.
\end{equation}

Recall $
\vnorm{(x,y)} = \enorm{\mathbf{A}^{-1}(x,y)}.
$
Since our eigenvectors are $\mathbf{v_j} = 
(\lambda_j, 1)$, 
\[
\mathbf{A} = 
\begin{bmatrix}
\lambda_1 & \lambda_2 \\ 1 & 1
\end{bmatrix},
\text{ \ and \ }
\mathbf{A}^{-1} = 
\frac{1}{\lambda_1 - \lambda_2}
\begin{bmatrix}
1 & -\lambda_2 \\ -1 & \lambda_1
\end{bmatrix}.
\]
First we show that we can set $C$ as claimed:
\begin{eqnarray*}
&& \abs{\lambda_1 - \lambda_2}^2 \vnorm{ \twovec{x}{y} }^2 
= \abs{\lambda_1 - \lambda_2}^2 \enorm{ \mathbf{A}^{-1} \twovec{x}{y} }
= \abs{x-\lambda_2 y}^2 + \abs{-x + \lambda_1 y}^2 \\
&\leq& 
2\abs{x}^2 +\abs{y}^2 ( \abs{\lambda_2}^2 + \abs{\lambda_1}^2) 
+2\abs{x} \abs{y}(\abs{\lambda_2} +\abs{\lambda_1}),
\text{ by triangle inequality},\\
 &\leq& 
2\abs{x}^2 +\abs{y}^2 ( \abs{\lambda_2}^2 + \abs{\lambda_1}^2) 
+(\abs{x}^2+\abs{y}^2)(\abs{\lambda_2} +\abs{\lambda_1}),
  \text{ by Equation~\ref{eqn:crossest}}, \\
  &=&   
\abs{x}^2 \( 2+\abs{\lambda_1}+\abs{\lambda_2} \)
+ \abs{y}^2 \( \abs{\lambda_2}^2+\abs{\lambda_1}^2 
              +\abs{\lambda_2}+\abs{\lambda_1} \), \\
 & \leq  &
 (\abs{x}^2 + \abs{y}^2) (2+\abs{\lambda_1}+\abs{\lambda_2}),
\text{ since } \abs{\lambda_1}, \abs{\lambda_2} \leq 1, \text{ so } \abs{\lambda_2}^2+\abs{\lambda_1}^2 \leq 2.
\end{eqnarray*}

Next we need to establish that $\enorm{\mathbf{u}} \leq D \vnorm{\mathbf{u}} = D \enorm{\mathbf{A}^{-1}\mathbf{u}}$,
for all $\mathbf{u} \in \Ct$. Since $\mathbf{A}$ is invertible, this is
is equivalent to:
$\enorm{\mathbf{A}\mathbf{u}} \leq D \vnorm{\mathbf{A}\mathbf{u}} = D \enorm{\mathbf{u}}$,
for all $\mathbf{u} \in \Ct$.  
The following establishes this latter statement with $D$ as claimed:
\begin{eqnarray*}
&& \enorm{\mathbf{A}\twovec{x}{y}}^2 = 
 \abs{\lambda_1 x + \lambda_2 y }^2 + \abs{x+y}^2, \\
& \leq &
\abs{x}^2 ( \abs{\lambda_1}^2 +1 ) +  \abs{y}^2 ( \abs{\lambda_2}^2 + 1) + 
2\abs{x} \abs{y} ( \abs{\lambda_1} \abs{\lambda_2} + 1 ) 
\ \text{(triangle ineq.)}, \\
& \leq &
\abs{x}^2 \( \abs{\lambda_1}^2 + \abs{\lambda_1} \abs{\lambda_2} +2 \)
+ \abs{y}^2 \( \abs{\lambda_2}^2 + \abs{\lambda_1} \abs{\lambda_2} +2 \)
\ \text{(Equation~\ref{eqn:crossest})}, \\
& \leq  &
\( \abs{x}^2 + \abs{y}^2 \) \( \lambda + \abs{\lambda_1} \abs{\lambda_2} +2 \), 
\text{ since we defined } \lambda = \max (\abs{\lambda_1}, \abs{\lambda_2}),\\
& \leq &
\( \abs{x}^2 + \abs{y}^2 \) \( \lambda + \abs{a} +2 \), 
\text{ since } \det (Df) = a, \text{ so } \abs{\lambda_1} \abs{\lambda_2} = a.
\end{eqnarray*}

\end{proof}

\begin{rem}
Since the eigenvectors are $(\lambda_j,1)$, the difference
$\abs{\lambda_1-\lambda_2}$ is the determinant of~$\mathbf{A}$.  This is
small if the angle difference between the eigenvectors is small. In 
this case, $C$ captures that the metric is skewed far from euclidean, so 
only a very small euclidean ball can fit inside a $\vsymb$-ball.  
Note also that $D$ is large only when the eigenvalues are large.  
Thus $D$ captures the strength of the contraction.  
\end{rem}

\subsection{Estimating the size of the sink basin}

Now that we know how to convert between the two norms, we are ready to 
take some measurements in the sink basin.  
To do so, we approximate $f$ by $L_p f$, the linearization of $f$ at $p=(z,z)$.  Recall
\[
L_p f \twovec{x}{y} = f\twovec{z}{z} + D_pf \twovec{x-z}{y-z}
= \twovec{z^2+c -2z(x-z) -ay}{x}.
\]

We first bound the error between $f$ and $L_p f$ in the 
$\vsymb$-norm.

\begin{lem} \label{lem:HLperror}
If  $\vnorm{\twovec{x-z}{y-z}} = r$, then  $\vnorm{f\twovec{x}{y} - L_p f\twovec{x}{y}}  \leq r^2 \( \frac{D^2}{C} \). $
\end{lem}

\begin{proof}
Let $(x,y) \in \Ct$ be such that
$ 
\vnorm{\twovec{x-z}{y-z}} = r,
$
for some $r>0$.

It is easy to compute the quadratic error in approximating $f$ with $L_p f$ 
in the euclidean metric, since
$
f\twovec{x}{y} - L_p f\twovec{x}{y} = \twovec{(x-z)^2}{0}.
$

We then convert to the $\vsymb$-norm, using 
Lemma~\ref{lem:vnormequiv} twice, to get:
\begin{eqnarray*}
& \vnorm{\twovec{(x-z)^2}{0}} 
 \leq  
 \frac{1}{C} \enorm{\twovec{(x-z)^2}{0}}
= 
\frac{1}{C} \abs{x-z}^2  
 \leq 
\frac{1}{C} \abs{x-z}^2 + \abs{y-z}^2 
& \\
& = 
\frac{1}{C} \enorm{\twovec{x-z}{y-z}}^2 
 \leq 
  \frac{D^2}{C} \vnorm{\twovec{x-z}{y-z}}^2 
= 
\(\frac{D^2}{C}\) r^2. &
\end{eqnarray*}
\end{proof}

Next, we show that in the $\vsymb$-norm, the linearization moves 
points closer to $p$ by a linear contraction. 

\begin{lem} \label{lem:Lperror}
If  $\vnorm{\twovec{x-z}{y-z}} = r,$   then 
$\vnorm{L_p f\twovec{x}{y} - \twovec{z}{z}}
 \leq \lambda r.$
\end{lem}

\begin{proof}
Since $p$ is fixed,
$
L_p f \twovec{x}{y} - \twovec{z}{z} = D_pf \twovec{x-z}{y-z}.
$
Now to work with $D_pf$, note that since the columns of $\mathbf{A}$ 
are the eigenvectors of~$D_pf$, we have 
$$
\mathbf{A}^{-1} D_pf \mathbf{A} = 
  \begin{bmatrix}   \lambda_1 & 0 \\ 0 & \lambda_2    \end{bmatrix}.
$$
Using this, and the fact that $\lambda = \max \( \abs{\lambda_1},\abs{\lambda_2} \)$, we get:
\begin{eqnarray*}
 \vnorm{D_pf\twovec{x-z}{y-z}} 
& = & \enorm{\mathbf{A}^{-1} D_pf \twovec{x-z}{y-z} }
 \\
& = & \enorm{ 
 \begin{bmatrix} \lambda_1 & 0 \\ 0 & \lambda_2 \end{bmatrix} 
    \mathbf{A}^{-1} \twovec{x-z}{y-z} }  \\
& \leq  & \lambda \enorm{\mathbf{A}^{-1} \twovec{x-z}{y-z}}
= \lambda \vnorm{\twovec{x-z}{y-z}} = \lambda r.
\end{eqnarray*}
\end{proof}

From the above lemmas and the triangle inequality, we immediately conclude:

\begin{lem} \label{lem:Hp}
If  $\vnorm{\twovec{x-z}{y-z}} = r,$  then 
$\vnorm{f\twovec{x}{y} - \twovec{z}{z}}
 \leq \lambda r + r^2 \( \frac{D^2}{C} \).$
\end{lem}

Now we can estimate the euclidean size of the sink basin.

\begin{prop} \label{prop:sinkbasin}
Let
\[
\tau = \frac{ \abs{\lambda_1 -\lambda_2}^2}
{(2+\abs{\lambda_1}+\abs{\lambda_2})(2+\lambda^2+\abs{a})}.
\]

Then the euclidean disk centered at $p$ of radius
$
r_p = \tau (1-\lambda)
$
is contained in the immediate sink basin of~$p$.
\end{prop}

\begin{proof}
Note $C^2/D^2 = \tau$.
We first show that the $\vsymb$-disk centered at $p$ of 
radius $s_p = (1-\lambda)(C/D^2)$, $\DD_{\vsymb}(p, s_p)$, is 
contained in the sink basin.  
For, if $\vnorm{(x-z,y-z)} =r \leq s_p$, then by Lemma~\ref{lem:Hp},
$
\vnorm{f\twovec{x}{y}-\twovec{z}{z}} \leq \lambda r + r^2 (D^2/C),
$
and by $r \leq s_p$ and definition of 
$s_p$ we get $\lambda r + r^2 (D^2/C) \leq r.$
Thus $f$ maps the disk $\DD_{\vsymb}(p, s_p)$ into itself, and 
every point in it closer to $p$ in the $\vsymb$-norm. Thus this 
disk is contained in the immediate sink basin.

Now if we use Lemma~\ref{lem:vnormequiv} to convert to a euclidean 
statement, we see that for $r_p := s_p \ C$ we have $r_p = s_p C = (1-\lambda)(C^2/D^2) = (1-\lambda)\tau,$ and 
$\DD_{e}(p,r_p) \subset \DD_{\vsymb}(p,s_p) \subset \{ 
\text{immediate sink basin} \}$.
\end{proof}

\subsection{Separating the \boxchtrcomps}

We now investigate the \boxchtr components, \textit{i.e.}, the strongly connected components of $\Gamma$, using their relation to the $\eps$-chain transitive components as given in
Corollary~\ref{cor:sccchaincomp}.  

First, given a sufficiently small constant $\xi$, we calculate an annular region
$\mathcal{A}_{\xi}$, contained in the immediate sink basin, in which the
contraction toward the fixed point causes iterates to move toward $p$ by a distance large enough to block $\xi$-chain recurrence, with respect to the $\vsymb$-norm.

\begin{lem} \label{lem:nonepsrec}
Let 
$
0 < \xi < (1-\lambda)^2 C /(4 {D^2}).
$
Then, in the $\vsymb$-norm,
 the $\xi$-chain transitive component that contains
the fixed sink is separated from the $\xi$-chain transitive component of any
other invariant set by a distance of 
$
(1-\lambda) C/D^2.
$
\end{lem}

\begin{proof}
Define $\mathcal{A}_{\xi}$ by
$
\mathcal{A}_{\xi} = \left\{ \twovec{x}{y} \colon 
   r_{-} < \vnorm{\twovec{x-z}{y-z}} < r_{+} \right\},
$
where $p=(z,z)$ is the fixed sink, and 
\[
r_{\pm} = \frac{C}{2D^2} \( (1-\lambda) \pm 
                        \sqrt{(1-\lambda)^2 - 4\xi D^2/C } \).
\]
We show that if $(x_0,y_0) \in \mathcal{A}_{\xi}$, it is not 
$\xi$-chain recurrent with respect to $\vnorm{\cdot}$.   Since $\mathcal{A}_{\xi}$ is of $\sigma$-width $(1-\lambda) C/D^2,$ this establishes the lemma.

Note that $r_{\pm}$ are the roots of the polynomial
\begin{equation} \label{eqn:qdef}
q(r) = (D^2/C) r^2 - (1-\lambda)r + \xi.
\end{equation}
Thus $\xi < (1-\lambda)^2 C / (4D^2)$ is precisely the condition that
needs to hold in order for the roots $r_{\pm}$ to be real, and thus positive.  
Hence,
\begin{equation} \label{eqn:qneg}
q(r) < 0, \text{ if } r \in (r_{-}, r_{+}).
\end{equation}

Now we show the contraction of $f$ in $\mathcal{A}_{\xi}$ is strong enough to block $\xi$-chain recurrence.  
Suppose $(x_0,y_0) \in
\mathcal{A}_{\xi},$ with $r = \vnorm{(x_0-z,y_0-z)}$.  
Let $n \in \ZZ^{+}$ and let
$\{ (x_1, y_1), \ldots, (x_{n-1},y_{n-1}) \}$ be any points such that 
$\vnorm{(x_{j},y_{j}) - f(x_{j-1},y_{j-1})} < \xi,$ 
for $0 < j \leq n-1$. 
We show that $\vnorm{(x_0,y_0) - f(x_{n-1},y_{n-1})} \geq \xi$ by
 first showing inductively that for $0 \leq j \leq n-1$, 
\begin{equation} \label{eqn:IHnonepsrec}
\vnorm{\twovec{x_j}{y_j}-\twovec{z}{z} } \leq r, 
\text{ so by Lemma~\ref{lem:Hp}, } 
\vnorm{f\twovec{x_j}{y_j} -\twovec{z}{z}} \leq \lambda r + r^2 
\frac{D^2}{C}.
\end{equation}

We have (\ref{eqn:IHnonepsrec}) for $j=0$ already.  Now let $0 < j \leq n-1$, 
and suppose we know (\ref{eqn:IHnonepsrec}) for $(x_{j-1},y_{j-1})$.  Then
\begin{eqnarray*}
&& \vnorm{\twovec{x_{j}}{y_{j}} - \twovec{z}{z} } 
\\
& \leq & \vnorm{ \twovec{x_{j}}{y_{j}} - f\twovec{x_{j-1}}{y_{j-1}} } 
        + \vnorm{f\twovec{x_{j-1}}{y_{j-1}} -\twovec{z}{z} } 
      \text{ (triangle ineq.),}  \\
& \leq & \xi + \lambda r + r^2 (D^2/C)  =  q(r) + r \leq r 
   \text{ (choice of } (x_{j}, y_{j}), \text{ Equations~\ref{eqn:qdef} and~\ref{eqn:qneg}).}
\end{eqnarray*}

Thus induction verifies (\ref{eqn:IHnonepsrec}).  In particular, (\ref{eqn:IHnonepsrec}) holds for $j=n-1$ and $j=0$.  But then
\begin{eqnarray*}
&& \vnorm{\twovec{x_0}{y_0} -f\twovec{x_{n-1}}{y_{n-1}} } \\
& \geq & \vnorm{ \twovec{x_0}{y_0} - \twovec{z}{z} } 
    - \vnorm{ f\twovec{x_{n-1}}{y_{n-1}} - \twovec{z}{z}} 
 \text{ (triangle ineq.)}
    \\
& \geq & r - \lambda r - r^2 (D^2/C)
 =  \xi - q(r) \geq \xi
\text{ (Equations~\ref{eqn:IHnonepsrec}, \ref{eqn:qdef}, and~\ref{eqn:qneg}). }
\end{eqnarray*}
Hence, $(x_0, y_0)$ is not $\xi$-chain recurrent in the $\sigma$-norm.
\end{proof}


A key component of the proof of Theorem~\ref{thm:babysinkJsep} is the following quantification of $\eta >0$ for which the $\eta$-chain recurrent set is {\em separating}, in a sense parallel to Definition~\ref{defn:sep}:  

\begin{defn}  \label{defn:chrecsep}
We call the $\eps$-chain recurrent set, $\chreceps$, {\em separating} if there are two chain transitive components, $\chainrec^j$ and $\chainrec^k$ (of $\chrec$), which lie in different $\eps$-chain transitive components of $\chreceps$.   In this case we say $\chreceps$ {\em separates } $\chainrec^j$ and $\chainrec^k$. 

Further, we call $\chreceps$ {\em fully separating} if it separates every pair of chain transitive components of $\chrec$.
\end{defn}

\begin{prop} \label{prop:babynonepsrec}
Suppose $f_{a,c}$ is a  \Henon mapping with an attracting fixed point $p$, with $\lambda_1 \neq \lambda_2$ eigenvalues of $D_pf$, and $\lambda = \max (\abs{\lambda_1}, \abs{\lambda_2})$. 
Let $\tau$ be as in Theorem~\ref{thm:babysinkJsep}.

If $\eta>0$ satisfies
$
\eta < \tau (1-\lambda)^2/4,
$
then $\chreceta$ is separating.  
In particular,  the $\eta$-chain transitive component containing the sink
is distinct from the $\eta$-chain transitive component of any
other invariant set.
\end{prop}

\begin{proof}
%
We have $\tau = C^2/D^2$.
We simply convert Lemma~\ref{lem:nonepsrec} to euclidean estimates,
using Lemma~\ref{lem:vnormequiv}.  Let $\eta = \xi C$, so that
\[
\eta = \xi C < \frac{(1-\lambda)^2}{4} \frac{C^2}{D^2}
= \frac{(1-\lambda)^2}{4} \tau,
\] 
as claimed in the statement of the proposition.

Define the set $\mathcal{S}_{\eta}$ to be simply the set
$\mathcal{A}_{\xi}$, from the
proof of Lemma~\ref{lem:nonepsrec}.  
Let $(x_0,y_0) \in \mathcal{S}_{\eta}$, and let
$n \in \NN$ and $\{ (x_1, y_1), \ldots, (x_{n-1},y_{n-1}) \}$ be any
points such that
$
\enorm{\twovec{x_{j}}{y_{j}} - f\twovec{x_{j-1}}{y_{j-1}}} < \eta, 
\ 0 < j \leq n-1.
$  
Then 
$
\vnorm{\twovec{x_{j}}{y_{j}} - f\twovec{x_{j-1}}{y_{j-1}}}
< \eta /C = \xi.
$
Thus as in the proof of Lemma~\ref{lem:nonepsrec}, we have
\[
\enorm{\twovec{x_0}{y_0} -f\twovec{x_{n-1}}{y_{n-1}}} 
\geq C \vnorm{\twovec{x_0}{y_0} -f\twovec{x_{n-1}}{y_{n-1}}} 
\geq C \xi = \eta.
\]
Thus, $(x_0,y_0)$ is not in $\chreceta$. Thus there exists a connected set $\mathcal{S}_{\eta}$ in the immediate sink basin which lies in $\Ct \setminus \chreceta$, hence $\chreceta$ separates the fixed sink from every other chain transitive component.
\end{proof}


Now we prove Theorem~\ref{thm:babysinkJsep}, by applying the estimates of Proposition~\ref{prop:babynonepsrec} to \boxchrec sets, to quantify which box size $\eps$ guarantees that the components of an $(\eps, \delta)$-\boxchrecmod $\Gamma$ separate the fixed sink from every other chain transitive component, in the sense of Definition~\ref{defn:sep}.  
%

Recall Theorem~\ref{thm:babysinkJsep} stated that 
for $M>1 $ such that $\delta < \eps/M$, and
$
\kappa 
:=
\left[ 1 + 1/M + \max \{1, (1-\lambda)\sqrt{\tau} +2\snorm{p} + \sonorm{a} 
\} \right], 
$
if
 $\eps < \frac{1}{2} \( -\kappa + \sqrt{\kappa^2+\tau(1-\lambda)^2 } \),$
then $\Gamma$ is separating.

\begin{proof}[Proof of Theorem~\ref{thm:babysinkJsep}]
Let $\boxcov = \boxcov({\Gamma})$ be the region in $\Ct$ covered by the vertices of the $(\eps, \delta)$-\boxchrecmod $\Gamma$.
First note by Proposition~\ref{prop:babynonepsrec}, if $\eta < \tau (1-\lambda)^2 /4$
there is a connected set 
$\mathcal{S}_{\eta}=\mathcal{A}_{\xi}$ in the 
immediate sink basin which lies in $\Ct \setminus \chreceta$.   
Now if we can calculate a bound on $\eps$ so that $\boxcov \cap  \mathcal{S}_{\eta}  \subset \chreceta$,
then we would have $\boxcov \cap \mathcal{S}_{\eta} = \emptyset$.  To maximize the bound on $\eps$, set $\eta = \tau (1-\lambda)^2 /4$.

Now Proposition~\ref{prop:sinkbasin} implies $\DD_{\vsymb}(p,s_p)$ is mapped into itself by $f$, thus
there are no edges in $\Gamma$ from boxes inside $\DD_{\vsymb}(p,s_p)$ to 
those outside of $\DD_{\vsymb}(p,s_p)$.  
Since $\mathcal{S}_{\eta} \subset \DD_{\vsymb}(p,s_p)$,
the \boxchtrcomp containing the sink must be distinct from the \boxchtrcomp containing 
$J$, whenever the box size $\eps$ is small enough that 
$\boxcov  \cap \mathcal{S}_{\eta}  \subset \chreceta$.

To get $\boxcov  \cap \mathcal{S}_{\eta}  \subset \chreceta$, first
recall that Theorem~\ref{thm:RinBoxCov}
calculates an $\eps'$ such that $\boxcov \subset \chrecepsp$.  
The theorem specifies that $\eps' = \delta + \eps(r+1)$, where $r$ is computed in 
Lemma~\ref{lem:imageboxsize} as $r=\eps+(2R'+\sonorm{a}),$ so $\eps' = \delta + \eps^2 + \eps(2R+\sonorm{a})+\eps$.
 Note for \Henon mappings, we do not need to require $\eps<1$, this was applied to higher order terms.  
 
By examining the proof of Lemma~\ref{lem:imageboxsize}, and restricting the estimates of that proof to apply only to $\boxcov \cap \mathcal{S}_{\eta}$, we find we can use a slightly better estimate for $r$.  Indeed, let $R^{+}$ be a bound on the box-norm, $\snorm{\cdot}$, of a 
point in $ \mathcal{S}_{\eta} $.  We calculate $R^{+}$ below. 
Also recall that by hypothesis we have $\delta < \eps/M$.
Set
\[
 \nu = \eps/M + \eps^2 
        + \eps \max \(1, 2R^{+} + \sonorm{a} \) + \eps.
\]
Then the proof of Lemma~\ref{lem:imageboxsize} implies  $\boxcov \cap 
 \mathcal{S}_{\eta}  \subset  \chrecnu$. 
If $\nu \leq \eta$, then $\boxcov \cap  \mathcal{S}_{\eta} \subset \chrecnu \subset \chreceta$.  So we just need $\eps$ small enough that $\nu \leq \eta$.

To compute $R^+$, we use that $ \mathcal{S}_{\eta} = \mathcal{A}_{\xi}$ in the 
$\vsymb$-norm is 
centered at the point $p$, 
and has outer radius 
\[
r_+
= \frac{C}{2D^2} \( (1-\lambda) +
                        \sqrt{(1-\lambda)^2 - 4\xi D^2/C } \).
\]
But at the maximum $\eta = \xi C$, the discriminant is zero, so we have
$r_+ \leq \frac{C}{2D^2} (1-\lambda)$.  Converting to the euclidean norm,
we have a bound of $Dr_+ \leq \frac{C}{2D} (1-\lambda) = \sqrt{\tau} 
(1-\lambda)/2$.  Now since the box-norm is less than euclidean, we get 
that if $(x,y) \in \mathcal{S}_{\eta}$, then $\snorm{(x,y)} \leq 
\sqrt{\tau} (1-\lambda)/2 + \snorm{p} =: R^{+}$.

Then $\nu = \eps^2 + \eps [ 1/M + \max\{1,( \sqrt{\tau} (1-\lambda) + 2\snorm{p}) + \sonorm{a}\}+1]$, hence notice $\nu = \eps^2 + \eps \kappa$.  Now to bound $\eps$ so that $\nu \leq \eta$, let $q(\eps) = \eps^2 + \eps \kappa - \eta$.
Then the roots of~$q$ are $(-\kappa \pm \sqrt{\kappa^2 +4\eta} )/2$.
Both roots are real, with one positive and one negative.  We seek $\eps 
>0$ small enough that $q(\eps) <0$, which is the same as $\eps$ smaller 
than the positive root $(-\kappa + \sqrt{\kappa^2 +4\eta} )/2$.  Since $\eta = \tau(1-\lambda^2)/4$, this is exactly the bound on $\eps$ claimed in the statement of the theorem.
\end{proof}

\begin{rem}  In the above proposition, note one does not need all boxes in $\boxcov$ to be of size 
$\eps$, but rather just the boxes in the immediate sink basin, computed in 
Proposition~\ref{prop:sinkbasin}.  Thus a selective subdivision procedure 
targeting the sink basin could be advantageous for speeding up separation.
\end{rem}

\subsection{One dimension}

All of the work of this section applies to $P_c(z)=z^2+c$
 in the case of a fixed sink $p$ with multiplier $\lambda = \abs{P_c'(p)}$.  In this case, 
we do not need the $\vsymb$-norm, so take $\tau=C=D=1$.  Then 
 the disk $\DD_e(p, (1-\lambda))$ is in the sink basin and 
for $\delta < (1-\lambda)^2/4$ the set
\[
\mathcal{A}_{\delta} = \{ z \colon r_{-} < \abs{z-p} < r_{+} \},
\text{ where }
r_{\pm} = \frac{1}{2} \( (1-\lambda) \pm \sqrt{(1-\lambda)^2-4\delta} \),
\]
is in $\CC \setminus \chrecdelta$.
In order to guarantee separation of~$J$ from the sink,
for
\[
\kappa = \(1+1/M+ (1-\lambda)+2\abs{p}) \) \text{ and }
\eta =(1-\lambda)^2/4, 
\]
we need boxes of side length $\eps$ satisfying
$
\eps < \(-\kappa + \sqrt{\kappa^2+4\eta}\) /2.
$

\subsection{The 3-1-map}

We now apply our estimates to the 3-1-map, to determine
how small the boxes would need to be to produce a \boxchrecset separating $J$ from
the fixed sink.  The results of this calculation are shown in Table~\ref{table:per31}.
\begin{table}
\begin{center}
\caption{Constants for sink/$J$ separation estimates for 
Example~\ref{exmp:per31}.}
\label{table:per31}
\begin{tabular}{|rcl|}
\hline
$p$ & $=$ & $(-0.612,-0.612)$ 
\\
$\lambda_1$ & $=$ & $-.885$ 
\\
$\lambda_2$ & $=$ & $-.34$ 
\\
$\lambda$ & $=$ & $.885$
\\ 
$\tau$ & $=$ & $0.029871571$
\\
$\tau (1-\lambda)$ & $=$ & $0.0034352307$ 
\\
$\kappa$ & $=$ & $2.5448759$
\\
$\eta$ & $=$ & $9.876288 \times 10^{-5}$
\\
$\eps$ & $<$ & $3.880793 \times 10^{-5}$ 
\\
\hline
\end{tabular}
\end{center}
\end{table}

Thus, a box side length less than $4 \times 10^{-5}$ would guarantee that the \boxchcn 
separate the fixed sink from $J$. 
%
%
However, this is several orders of magnitude
smaller than the best we could compute with current resources ($0.015$).  
Also, note that the guaranteed euclidean disk contained in the sink basin
is only of radius $0.0034$.  But visually inspecting this example suggests the immediate basin is much larger.  
This suggests a computer should be able to rigorously separate the dynamics of $J$ from the sinks, but more sophisticated techniques are needed.

\section{Examples}
\label{sec:examples}

In this section we present several examples of applying the \boxchcn to 
\Henon mappings. Every example uses the procedure of eliminating $\boxcov_0$-escaping boxes, and some examples use \inKsel subdivision.  We also present an example of the \boxchcn applied to a polynomial map of $\CC$.

The computations described in this section were run on a Sun Enterprise E3500
server with $4$ GB of RAM and four processors, each $400$MHz UltraSPARC (though the multiprocessor was not used).  \begin{footnote}{The server was obtained by the Cornell University Mathematics Department through an NSF SCREMS grant.}\end{footnote}
When computations became overwhelming, memory usage was the limiting 
factor.

To measure the accuracy of an $(\eps, \delta)$-\boxchrecmod of $f$, $\Gamma$, we compute the bounds $\eps'$ and $\delta'$ given in Theorem~\ref{thm:babyRinBoxCov} such that $\chrecdeltap \subset \boxcov({\Gamma}) \subset \chrecepsp$.  As we discussed in Corollary~\ref{cor:RinBoxCov}, the maximum error $s(F)$ between the ideal function $f$ and an implemented interval extension of it, $F$,
 evaluated for any box in the model, must be added to $\eps'$ to get the corresponding result for the actually computed model. However since at worst 
$s(F)$ is still much less than $\eps'$, we neglect discussing this factor in the examples presented below.

Table~\ref{table:examples} contains more detailed data for all of the examples of this section.

\subsection{Polynomial maps of $\CC$}

Recall that all of our results can be reformulated for polynomial maps of $\CC$ of degree $d>1$.
For example, consider the cubic polynomials $P_{a,c}(z) = z^3-3a^2 z+c$.  One can check that $R'=2$
suffices for $\abs{c} < 2$ and $\abs{a} \leq \sqrt{2/3}$. 
Below we describe an example of an interesting \boxchrecmod of a cubic polynomial, computed with Hypatia.

\begin{exmp} \label{exmp:cubdoublebas}
The cubic polynomial $P_{a,c} (z)= z^3-3a^2 z+c,$ with $c=-.19+1.1i, a=0.1i$, has 
an $\chainrec$ consisting of $J$ and a 4-cycle. Figure~\ref{fig:cubdoublebas} 
\begin{figure}
\begin{center}
\drawfigcubdoublebas
\caption
{A \boxchrecset for  $P_{a,c} (z)= z^3-3a^2 z+c,$ with $c=-.19+1.1i, a=0.1i$, and boxes from a $2^{10} \times 2^{10}$  grid on $[-2,2]^2$. 
Here $\chainrec$ is $J$ and a 4-cycle.  The black regions form the \boxchtrcomp of the 4-cycle. 
The origin lies in the center of the largest black region (in the center of the figure).  The \boxchtrcomp containing $J$ is shaded two tones, to heuristically  illustrate J.  
The distance from the origin to the end of any of the three limbs of $J$ is approximately $1.3$. 
The small shaded islands spread between $J$ and the 4-cycle neighborhood are
several \boxchtrcomps which do not intersect $\chainrec$, but are contained in some $\chrecepsp$. 
}
\label{fig:cubdoublebas}
\end{center}
\end{figure}
shows the \boxchrecset which we computed using Hypatia. 
To get this approximation, we uniformly subdivided to obtain boxes from a $2^{10} \times 2^{10}$  grid on $[-2,2]^2$.  Note there are several small \boxchtrcomps spread between $J$ and the 4-cycle, which cannot contain any points of $\chainrec$, but by Corollary~\ref{cor:sccchaincomp} these \boxchtrcomps are contained in some $\chrecepsp$.  
\end{exmp}

\subsection{Drawing pictures of \boxchrecsets for \Henon mappings}
\label{sec:drawboxcovwup}
Since applying the \boxchcn is an iterative process, we need feedback after each step $(n)$ to help us determine when to halt.  In this paper, we consider a sucessful construction to be one which arrives at a fully separating \boxchrecmod  (recall Definition~\ref{defn:sep}).
 
To check heuristically whether a \boxchrecmod $\Gamma$ is separating, we take advantage of the fact that the unstable manifold, $W^u(p)$, of a saddle fixed point, $p$, can be (approximately) explicitly parameterized by the plane (details are given in Section~\ref{sec:henwuppics}), and we
sketch the slice of $\boxcov$ lying in some $W^u(p)$.  

This parameterization identifies $p$ with the origin, and conjugates the \Henon map to multiplication by the unstable eigenvalue. Hence a  parameterized picture of $J \cap W^u(p)$ is invariant up to scaling by this eigenvalue, so the restriction to any neighborhood of $p$, $\snbd{p}{r}$, shows a fundamental domain.  A \boxchrecset $\boxcov$ contains boxes of a certain size, so a parameterized $\boxcov \cap W^u(p)$  is not quite invariant under scaling, but for many choices of $r$, a sketch of $\boxcov \cap W^u(p) \cap \snbd{p}{r}$ will show us whether the box chain transitive components have separated $J$ from the sinks, and give an idea of how well $\boxcov$ is approximating $J$, especially compared to a picture generated the same way with larger boxes (on a previous step).

To determine the coloration
of a pixel, we check whether the pixel intersects some boxes in $\boxcov$.  
The simplest picture would entail coloring a pixel black if it hits $\boxcov$, or white if not.
Going one step further, we illustrate the multiple \boxchtrcomps of $\Gamma$ by using multiple shades of grey.  However,
since our picture is a parametrization of a one complex
dimensional manifold which does not line up with the axes in $\Ct$, a
pixel may hit more than one box, and in more than one \boxchtr component. So, we use a color palette which has different shades of grey for each \boxchtr component, and a distinctive shade (like black) if the pixel hits more than one \boxchtr component. 

We can also enhance our sketch in order to see how well a given \boxchrecset approximates $\chainrec$.
Theorem~\ref{thm:JinR} (Section~\ref{sec:dpspaces}) states that $J \subset \chainrec \subset K$,
and  if $W^u(p)$ is the unstable manifold of any saddle periodic point, $p$, then clearly $W^u(p) \subset K^{-}$.
Hence we slightly lighten the pixels which seem to be
in $K^+$; \textit{i.e.}, the center point of the pixel stays small
after several iterates of the map.
In this way we can check visually how well $\boxcov$ approximates~$J$.  

Note that since $W^u(p) \subset K^-$, the attracting periodic orbits of the map do not lie in $W^u(p)$. However, their basins of attraction can intersect $W^u(p)$.  Thus in sketching \boxchrecsets in this manner, we usually see part of the \boxchtrcomps corresponding to the sink orbits. 

\subsection{Complex \Henon mappings}

Recall from Section~\ref{sec:intro} that if $f_{a,c}$ is a \Henon mapping with
$a$ sufficiently small and $c$ is such that the polynomial $P_c(z) = z^2+c$ is hyperbolic (thus stable under perturbation), then $f|_J$ is topologically conjugate to the
function induced by $P$ on the inverse limit
${\lim_{\leftarrow} (J, P)}$ (\cite{HOV2}).  In this case, we say that  {\em $f$ is described by $P$}, or simply that {\em $f$ exhibits one dimensional behavior}.

\begin{exmp}[The alternate basilica] \label{exmp:altpertwo}
Recall from Example~\ref{exmp:introaltpertwo} 
 that
the \Henon mapping $f_{a,c}$ with $c=-1.1875, a=.15$,
seems to be one of the simplest \Henon mappings
which does not exhibit  one dimensional behavior.
Computer evidence suggests that $\chrec$ consists of $J$ and an attracting two-cycle.
We found the most efficient \boxchcn for separating the sink from $J$ 
was to subdivide $\boxcov_0$ initially with a $(2^6)^4$ grid, then perform \inKsel subdivision.  In this way we acheived the left side of Figure~\ref{fig:altpertwo} (shown in Section~\ref{sec:intro}).  To improve the approximation, we performed  \inKsel subdivision once more, and obtained the right side of Figure~\ref{fig:altpertwo}  (shown in Section~\ref{sec:intro}). In these figures, note the \boxchtrcomps skirting the inner edge of the more refined $\boxcov'$. These would not be present for smaller box size, \textit{i.e.}, these components do not contain any points of $\chainrec$, though they are contained in some~$\chrecepsp$.
\end{exmp}


Next we examine a different type of diffeomorphism: a \textit{horseshoe}. 

\begin{defn} \label{defn:horse}
A {\em horseshoe} is a diffeomorphism $f$ such that $f|_J$
is topologically conjugate to the left shift operator on~$\Sigma_2$ (the
symbol space of bi-infinite sequences of 0's and 1's).
The {\em horseshoe locus} in the \Henon parameter space, $\mathcal{H}$,
is the hyperbolic component of parameter space containing the set
of horseshoes.
\end{defn}

Since the full left shift on~$\Sigma_2$ is the inverse limit of the 
one-sided left shift, horseshoes exhibit one dimensional 
behavior.   All horseshoes are
topologically the same as Smale's horseshoe, \textit{i.e.}, the Julia set
is a Cantor set, and the dynamics are well-understood.

However, pictures of the \Henon parameter space produced by SaddleDrop (\cite{HubKarl})
suggest that the topology of the horseshoe locus is quite complicated.  There are intriguing conjections about the horseshoe locus, which motivate the study of complex \Henon horseshoes.  Below we describe one example of applying the \boxchcn to a horseshoe.

\begin{exmp} \label{exmp:complexhorse}
The \Henon mapping with $c=-2.75, a=-.74$ appears to be a \textit{complex horseshoe}---though parameters $c,a$ are real, the horseshoe for the map of $\Ct$ is not contained in $\Rt$.
For horseshoe diffeomorphisms, there is no sink.  Thus there is no clear option 
for selective subdivision. 
We simply had
Hypatia uniformly subdivide the boxes, to obtain a \boxchrecmod consisting of boxes from a $(2^6)^2$ grid on $\boxcov_0=\snbd{0}{2.84}$.  
Since $J$ is a Cantor set, it is easiest to visualize $J$ via the FractalAsm picture shown on the left side of Figure~\ref{fig:complexhorse}. 
\begin{figure}
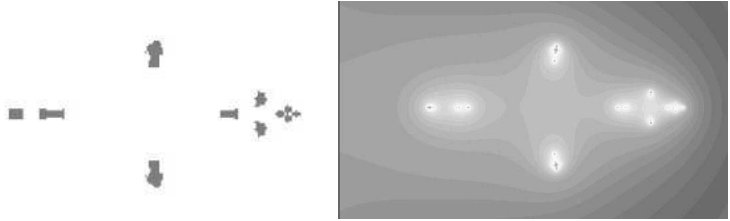

\drawfigcomplexhorse  \drawfigcomplexhorseFA
\caption{In these pictures, $f_{a,c}$ is the \Henon mapping with $a=-.74, c=-2.75$. This map appears to be a horseshoe, thus $\chainrec=J$ appears to be a 
Cantor set. On the right is a parametrization of $W^u(p)$, where black is $W^u(p) \cap K^+$.   On the left  is a \boxchrecset with boxes of size $2R'/2^{6}$, for $R'=2.84$.}
\label{fig:complexhorse}
\end{figure}
 The \boxchrecset is shown on the right side of 
Figure~\ref{fig:complexhorse}. 
Actually, this cantor set was sparse enough that we were able to refine uniformly to a $(2^{10})^2$ grid, but at this point the picture was nearly impossible to see, since the set is very small.
\end{exmp}

\subsection{Real \Henon mappings}

The \Henon mapping is widely studied as a diffeomorphism of $\Rt$, with $a,c$ real parameters. 
In fact, for some complex \Henon mappings with real parameters, $\chrec$ lies in $\Rt$, and the complex dynamics can be completely described by studying the restriction to $\Rt$ (see \cite{BS1}).  
Observe that all of the results of this paper apply immediately to the real setting.  We implemented the \boxchcn for real \Henon mappings in Hypatia; below is one example of a \boxchrecset for a real \Henon mapping.

\begin{figure}
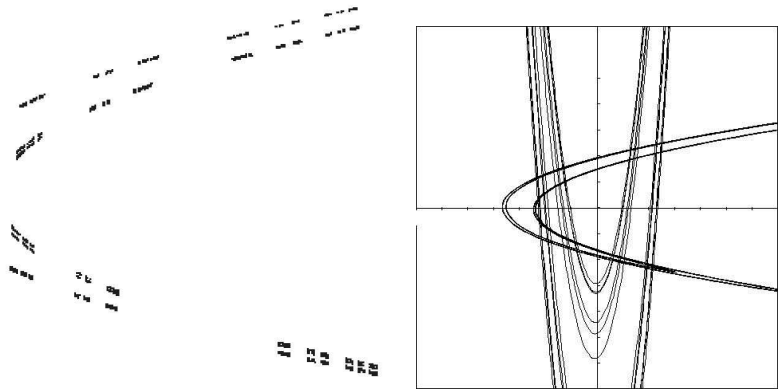

\begin{center}
\drawfigHrealhorse 
\ \ \  
\drawfigHrealhorseW
\caption{The \Henon mapping $f_{a,c}$ with $(a,c)=(-.25, -3)$ appears to be a real horseshoe, thus $\chainrec = J \subset \Rt$ is a Cantor set.  
On the right is a sketch of the stable and unstable manifolds of a saddle fixed point, $W^s(p) \cup W^u(p)$, shown in the region $[-7,7]^2$. On the left is a \boxchrecset $\boxcov$ with boxes of size $2R'/2^{7}$, where $R' = 2.56$. In fact, $\boxcov$ is a neighborhood of the intersection $W^s(p) \cap W^u(p)$.}
\label{fig:Hrealhorse}
\end{center}
\end{figure}

\begin{exmp} \label{exmp:Hrealhorse}
The \Henon mapping $f_{a,c}$ with $(a,c)=(-.25, -3)$ appears to be a \textit{real
horseshoe}, in that $\chainrec = J$ seems to lie in $\Rt$.  Thus a parameterized unstable manifold picture of the Julia set would simply show a Cantor set lying in the real axis.
With Hypatia we constructed a \boxchrecset using boxes from a $(2^7 \times 2^7)$ grid on $\snbd{0}{2.56} \subset \Rt. $ 
This \boxchrecset (shown in $\Rt$) is on the left in Figure~\ref{fig:Hrealhorse}.
\end{exmp}

The reader familiar with the study of real \Henon mappings may recognize that the \boxchrecset (on the left of Figure~\ref{fig:Hrealhorse}) appears to show an approximation to the intersection of the stable and unstable manifolds, $W^u(p)$  and $W^s(p)$ for some saddle periodic point $p$.  Indeed, the right side of Figure~\ref{fig:Hrealhorse} shows a sketch of such manifolds.  This observation holds because for any \Henon mapping, $\overline{W^u(p)} \subset J^{-}$ and $\overline{W^s(p)} \subset J^{+}$ (see Section~\ref{sec:bckgnd}).  Thus $J = J^+ \cap J^- = \overline{W^s(p)} \cap \overline{W^u(p)}$.

\begin{sidewaystable}
\caption{Data for the \boxchmods constructed in the Examples of Section~\ref{sec:examples}. Here $\Gamma'$ denotes the \boxchtrcomp of $\Gamma$ which contains $J$, and the \textit{box grid depth} for a box is the number $n$ such that the box is of size $2R'/2^n$.  If a \boxchmod contains boxes of multiple sizes, then multiple box grid depths are listed. 
 Recall Example~\ref{exmp:cubdoublebas} is of a cubic polynomial $P_{a,c}(z) = z^3-3a^2 z+c$, and the other examples are of \Henon mappings $f_{a,c}(x,y) = (x^2+c-ay,x)$.}
\label{table:examples}
\begin{center}
\begin{tabular} {| |    l  l  ||c||c| c |c|c||}
\hline \hline
Example & & \ref{exmp:cubdoublebas} & 
\multicolumn{2}{c|}{\ref{exmp:altpertwo}} & \ref{exmp:complexhorse} & \ref{exmp:Hrealhorse} \\ \hline

Figure & & \ref{fig:cubdoublebas} & 
\ref{fig:altpertwo} & \ref{fig:altpertwo} & \ref{fig:complexhorse} & \ref{fig:Hrealhorse}
\\ \hline

params. & $c$        & $-.19+1.1i$  &     
 \multicolumn{2}{c|}{$-1.1875$}   &   $-2.75$    &     $-3$                 \\ \cline{2-7}
	      & $a$        & $0.1i$           &
	      \multicolumn{2}{c|}{$ 0.15$}       &   $-.74$     &    $-.25$             \\ \hline
sink period &        &  $4$   &  
\multicolumn{2}{c|}{$2$}  & N/A  &  N/A     \\ \hline
$\eps'$ s.t. & $\boxcov \subset \chainrec({\eps'+s(F)})$ 
     & $.028$ & 
     $0.30$ & $0.30$ & $0.68$ & $0.26$ \\ \hline
$\delta'$s.t. & $\chrecdeltap \subset \boxcov$ 
      & $1.5 \times 10^{-6}$ & 
      $6\times 10^{-6}$ & $3\times 10^{-6}$ & $1.2 \times 10^{-5}$ & $6.3\times 10^{-6}$  \\ \hline
$R'$    &   &  $2.1$   &  
\multicolumn{2}{c|}{$1.9$} & $2.84$  & $2.57$ \\ \hline
box grid depth, $n$  & & $11$  & 
$6,7$  & $6, 7,8$  & $6$ & $7$\\ \hline
 box size & &  $0.002$   & 
 $.06, .03$ &$ .06, .03,.015$ & $0.09$ & $0.04$ \\  \hline
\# $\Upsilon$ boxes  & original  & $98$ & 
$184$&  $682$ & $60$ & $20$\\ \cline{2-7}
($1000$s) &  $\boxcov$-escaping &  $0$ & 
$116$ & $417$ & $43$ &   $15$ \\ \hline
$\Upsilon$ size & boxes & $98$  & 
$68$ & $265$ & $17$ & $4$\\ \cline{2-7}
($1000$s)          & edges &  $1{,}300$  &  
$2{,}500$ &$12{,}400$ & $830$ &  $120$\\ \hline
$\Gamma'$ size &   boxes &  $60$   & 
$53$& $182$ & $10$ &   $2.4$  \\ \cline{2-7}
 ($1000$s)          &  edges &  $780$   & 
 $2{,}500$ & $7{,}800$ & $500$ &  $75$\\ \hline
runtime (min.)   & &    $<5$    & 
& $20$ &  $<1 $& $<1$\\  \hline
RAM (MB)          & & $<200$  & 
&$900$ & $20$ & $20$ \\ \hline

\hline
\hline
\end{tabular}
\end{center}
\end{sidewaystable}

\appendix

\section{Background on \Henon mappings}
\label{sec:bckgnd}

\subsection{The \Henon Family}
\label{sec:henonintro}

Polynomial diffeomorphisms of $\Ct$ necessarily have polynomial inverses,
thus are often called polynomial automorphisms.
Friedland and Milnor (\cite{FM}) showed that
polynomial automorphisms of $\Ct$ break down into two categories.
\textit{Elementary} automorphisms have simple dynamics, and are
polynomially conjugate to a diffeomorphism of the form $(x,y) \mapsto
(ax+b, cy+p(x))$ ($p$ polynomial, $a,c \neq 0$). \textit{Nonelementary}
automorphisms are all conjugate to
finite compositions of \textit{generalized \Henon
mappings}, which are of the form $f(x,y) = (p(x)-ay,x)$, where $p(x)$ is a
monic polynomial of degree $d>1$ and $a \neq 0$.

To clarify the situation, one can define a \textit{dynamical degree} of a 
polynomial automorphism of $\Ct$.  If \textit{deg}$(f)$ is the maximum of the 
degrees of the coordinate functions, the dynamical degree is 
\[ 
d = d(f) = \lim_{n \to \infty} (\textit{deg}(f^{n}))^{1/n}. 
\] 
This degree is a conjugacy invariant.  Elementary automorphisms have 
dynamical degree $d=1$.  A nonelementary automorphism is conjugate to some
automorphism whose polynomial degree is equal to its dynamical degree. 
Without loss of generality, we assume such $f$ are finite
compositions of generalized \Henon mappings, rather
than merely conjugate to mappings of this form.

Thus, the quadratic, complex \Henon family $f_{a,c} (x,y)=(x^2+c-ay,x)$
represents the dynamical behavior of the simplest class of nonelementary
polynomial automorphisms; those of dynamical degree two.  We shall state
results assuming $f$ is a polynomial diffeomorphism of $\Ct$ with
$d(f)>1$, and often concentrate on the illustrative example of $f_{a,c}
(x,y)=(x^2+c-ay,x)$.

\subsection{Invariant sets of interest}
\label{sec:dpspaces}

The {\em chain recurrent set}, $\chainrec$, and
the {\em Julia set}, $J$, 
are both attempts at
locating the points with dynamically interesting behavior.
$\chainrec$ can also be decomposed into components which do not
interact with one another. 

First we recall the key concepts associated to chain recurrence.
An {\em $\eps$-chain} of length $n>1$ from $p$ to $q$ is a 
  sequence of points
  $\{p=x_1, \ldots , x_n=q\}$ such that  $\abs{f(x_k) - 
  x_{k+1}} < \epsilon$ for $1 \leq i \leq n-1.$
A point $p$ belongs to the {\em $\eps$-chain recurrent set}, 
  $\chreceps$, of a function~$f$ if there is an
  $\eps$-chain from $p$ to $p$.
The {\em chain recurrent set} is $\chainrec = \cap_{\eps>0} 
\chreceps.$
A point $q$ is in the {\em forward chain limit set of a point $p$},
$\chainrec(p)$, if for all $\eps >0$, for all 
$n\geq 1$, there is an $\eps$-chain from $p$ to $q$ of length 
greater than $n$.
Put an equivalence relation on $\chainrec$ by: 
$p \sim q$ if $p \in \chainrec(q)$ and $q\in\chainrec(p)$.
Equivalence classes are called {\em chain transitive components}.
We can define $\chreceps(p)$ and $\eps$-chain transitive
components analogously. For ease of notation, we will sometimes refer
to $\chainrec$ as $\chreceps$ for $\eps=0$, or $\chainrec(0)$.
Note that $\chainrec$ is closed and invariant, and if $\eps < \eps'$,
then $\chainrec \subset \chreceps
 \subset \chrecepsp$.  

Chain recurrence is quite natural to rigorously study using a computer.  A \boxchrecmod $\Gamma$ is an approximation to the dynamics of $f$ on $\chainrec$, 
and the connected components of $\Gamma$ are approximations to the chain transitive components.  
This is made precise in Section~\ref{sec:boxcover}.

For a polynomial map $f$ of $\CC$, the {\em filled Julia 
set}, $K$, is the set of points
whose orbits are bounded under $f$;  the {\em Julia set},
$J$,  is the topological boundary of $K$. 
For a polynomial diffeomorphism $f$ of $\Ct$, there are 
corresponding Julia sets:
$K^+ (K^-)$ is the set of points whose orbits are bounded under
$f (f^{-1})$ 
and $K = K^+ \cap K^-$ is called the \textit{filled Julia set};
$J^{\pm} = \partial K^{\pm}$ (the topological boundary) 
and $J = J^+ \cap J^-$ is called the \textit{Julia set}.
The Julia set can be easily sketched by computer, and is also extensively used in formulating theoretical results for complex \Henon mappings.

%

Bedford and Smillie show the following relationships between  $\chainrec$ and $J$

\begin{thm} [\cite{BS2}] \label{thm:JinR}
Let $f$ be a polynomial diffeomorphism of $\Ct$, with $d(f)>1$. 
\begin{enumerate}
\item Then $J \subset
\chainrec \subset K$ and $J$ is contained in
a single chain transitive component of $\chainrec$.
\item Assume further that
$\abs{\text{det}Df} < 1$.  Let $O_j$ for $j=1, 2, \ldots$ denote
the sink orbits of $f$.
\begin{enumerate}
\item Then $\chainrec$ is the set of bounded orbits (in forward/backward
time) not in punctured basins, where
if $p$ is a sink, the {\em punctured basin} of $p$
is $W^s(p) - p$, and
\item the chain transitive components are the sink orbits, $O_j$,
and the set $\chainrec - \cup_j O_j$.
\end{enumerate}
\end{enumerate}
\end{thm}

\subsection{Drawing Meaningful pictures for maps of $\Ct$}
\label{sec:henwuppics}

Filled Julia sets are the invariant sets which can be easily sketched by
computer, on any two-dimensional slice.  Hubbard has suggested the
following method for drawing a dynamically significant slice of the Julia
set of a \Henon mapping, by parameterizing an unstable manifold.  This method has been
implemented by Karl Papadantonakis into a program called FractalAsm, available for download at \cite{CUweb}.
We also use this method to draw pictures of collections of boxes, see Sections~\ref{sec:intro} and~\ref{sec:examples} for examples.

Let $f$ be a diffeomorphism of $\Ct$.  If $p$ is a periodic point of 
period $m$, and the eigenvalues
$\lambda, \mu$ of $D_pf^m$ satisfy $\abs{\lambda} > 1 > \abs{\mu}$ (or
vice-versa), then $p$ is a {\em saddle periodic point}.
The large (small) eigenvalue 
is called the unstable
(stable) eigenvalue.
If $p$ is a saddle periodic point,
then the {\em stable manifold of $p$} is
$ 
W^s(p) = \{q \colon d(f^n(q),f^n(p)) \rightarrow 0 \text{ as }
n\rightarrow \infty\}, 
$
and the {\em unstable manifold of $p$} is
$ 
W^u(p) = \{q \colon d(f^{-n}(q),f^{-n}(p)) \rightarrow 0 \text{ as }
n\rightarrow \infty\}. 
$
If $p$ a saddle periodic point of $f$, then $W^u(p) \ (W^s(p))$ is biholomorphically 
equivalent 
to $\CC$, and on $W^u(p) \ (W^s(p))$, $f$ is conjugate 
to multiplication by the unstable
(stable) eigenvalue of $D_pf$.

Let $f$ be a \Henon mapping.
When $\abs{a} \neq 1$, except on the curve of equation $4c = (1+a)^2$,
the map $f_{a,c}$ has at least one saddle fixed point, $p$, 
(\cite{HubKarl}). The unstable manifold 
$W^u(p)$ has a natural parametrization $\gamma \colon \CC \to W^u(p)$ 
given by
\[
\gamma(z) = 
\lim_{m\to\infty} \gamma_m(z) = 
\lim_{m\to\infty} 
f^m \( p + \frac{z}{\lambda_1^m}\mathbf{v_1} \),
\]
where $\lambda_1$ is the unstable eigenvalue of $D_pf$ and
$\mathbf{v_1}$ is the associated eigenvector.
This parametrization has the property that
$
f(\gamma(z)) = \gamma(\lambda z),
$
and any two parametrizations with this property differ by scaling the 
argument.

To parameterize $W^u(p)$, we approximate $\gamma$ by some $g=\gamma_N$ in
a region in the plane: $B=\{z=x+iy \colon a \leq x \leq b, c \leq y \leq
d\}$.  
Observe that since $W^u(p) \subset K^-$, to sketch $K$ in
$W^u(p)$, we need only sketch $K^+$.
Thus for each pixel $Z\in B$, if $f^{n}(g(Z))$ is bounded by some $R$ for all $n<N$, we
guess $g(Z)\in K^+$ and color $Z$ black.  Otherwise, color according to which
iterate $f^n(g(Z))$ first surpassed~$R$.

\section{Rigorous Arithmetic}
\label{sec:IA}

On a computer, we cannot work with real numbers; instead we work over
the finite space $\FF$ of numbers representable by binary floating point
numbers no longer than a certain length.  
For example, since the number 
$0.1$ is not a dyadic rational, it has an infinite binary expansion.  Thus 
the computer cannot encode $0.1$ exactly. 
Interval arithmetic (IA) provides a method for maintaining rigor in computations, and also is natural and efficient for manipulating boxes.  The basic objects of IA
 are closed intervals, $[a] = [\underline{a}, \bar{a}] \in \II\KK$, with end points in some fixed
 field, $\KK$.   An arithmetical operation on two intervals produces a resulting interval
which contains the real answer.  For example,
\begin{align*}
[a] + [b] & := 
\left[ \underline{a}+\underline{b}, \bar{a}+\bar{b} \right] \\
[a] - [b] & := 
\left[ \underline{a}-\bar{b}, \bar{a}-\underline{b} \right] 
\end{align*}
Multiplication and division can also be defined in IA.

Since an arithmetical operation on two computer numbers in $\FF$ may not have a 
result in $\FF$,  in order to implement rigorous IA we must round
outward the result of any interval
arithmetic operation, \textit{e.g.} for $[a],[b] \in \II\FF$,
$$
[a] + [b]  := 
\left[ \left\downarrow \underline{a}+\underline{b} \right\downarrow, 
         \left\uparrow  \bar{a}+\bar{b} \right\uparrow \right],
$$
where $\left\downarrow x \right\downarrow$ denotes the largest number in $\FF$ 
that is 
strictly less than $x$ (\textit{i.e.}, $x$ rounded down), and 
 $\left\uparrow x 
\right\uparrow$ 
denotes the smallest number in $\FF$ that is strictly greater than $x$ 
(\textit{i.e.}, $x$ 
rounded up).  This is called IA with \textit{directed rounding}.

For any $x\in \RR$, let Hull$(x)$ be the smallest interval in $\FF$ which contains $x$.
That is, if $x \in \FF$, then Hull$(x)$ denotes $[x, x]$.  If $x \in \RR \setminus \FF$, then 
 Hull$(x)$ denotes $\left[ \left\downarrow x \right\downarrow, 
\left\uparrow x \right\uparrow \right]$.
Similarly, for a set $S\subset \RR$, we say Hull$(S)$ for the smallest  interval containg $S$. Whether Hull$(S)$ is in $\IR$ or $\IF$ should be clear from context.

In higher dimensions, IA operations can be carried out component-wise, on 
{\em interval vectors}.  
So if $x\in \Rn$, then Hull$(x) = \text{ Hull}(x_1) \times \cdots \times
\text{ Hull}(x_n)$, and if $S \subset \Rn$, then $\text{Hull}(S)$ is the smallest vector in $\IF^n$ (or $\IR^n$) containing $S$.
Note to deal with intervals in $\Cn$ we simply identify $\Cn$ with $\RR^{2n}$. 
Thus a box in $\CC = \Rt$ is an interval vector of length two.
Our extensive use of boxes is designed to 
make IA calculations natural.

To compute the image of a point $x$ under a map $f$ using IA, first convert $x$ to an interval vector $X = \text{ Hull}(x)$, then use a (carefully chosen) combination of the basic arithmetical operations to compute an interval vector $F(X)$, such that we are guaranteed that Hull$(f(x)) \subset F(X)$.  

\begin{defn} \label{defn:intext}
Let $f \colon \Rn \to \Rn$ be continuous.  An \textit{interval extension} of $f$, $F = F(f)$, is a function which maps a box $B$ in $\Rn$ to a box $F(B)$ containing $f(B)$, \textit{i.e.},  $F(B)\supset \text{Hull}(f(B))$.
\end{defn}
Usually, 
we would like $F(B)$ to be as close as possible to Hull($f(B)$).
We shall not discuss how to find the best $F$. However, in this paper we do distinguish between $F(B)$ and Hull$(f(B))$ when stating our results, to make clear the difference between theory and practice.

Each time an arithmetical calculation is performed, one must think carefully
about how to use IA. For example, IA is not distributive.  Also, it can easily create large error propagation.
For example, iterating a polynomial map or diffeomorphism on an interval vector (which is not very close to an 
attracting period cycle) 
will usually produce a very large interval vector after only a few iterates.  
That is, if 
$B$ is a box in $\Rt = \CC$, and one attempts to 
compute a box containing $f_c^{10}(B)$, for $f_c(z)=z^2+c$, by:
\begin{tabbing}
\hspace*{.25in} \= \hspace{.25in} \= \hspace{.25in} \kill
\>for $j$ from $1$ to $10$ do \\
\>\>$B = F_c(B)$
\end{tabbing}
then the box $B$ will likely grow so large that its defining bounds become 
machine $\infty$, \textit{i.e.}, the largest floating point in $\FF$.
Similarly, if $f$ is a \Henon mapping, one would also never want to try to compute
$
D_{B_n}f \circ \cdots \circ D_{B_1}f \circ D_{B_0}f (\mathbf{u}),
$
for a vector $\mathbf{u} \in \Ct$, since the entries would blow up.

We use IA for all the rigorous computations in the computer program Hypatia.  The
IA routines were all provided by the PROFIL/BIAS package, available at 
\cite{PBIA}.  For further background on interval arithmetic, see  \cite{ GenIA, MooreIA1, MooreIA2}.


\bibliographystyle{plain}
\bibliography{SLHChainRec}



\end{document}